\documentclass[11pt,a4paper,reqno]{amsart}

\usepackage[dvipsnames]{xcolor}
\usepackage[color,all]{xy}
\usepackage{amsmath,amsfonts,amssymb,amsthm}
\usepackage{tikz}
\usetikzlibrary{matrix,positioning,decorations.markings,arrows,decorations.pathmorphing,backgrounds,fit,positioning,shapes.symbols,chains,shadings,fadings,calc}
\usepackage{tikz-cd}
\usepackage[colorlinks,linkcolor=blue!72,anchorcolor=orange,citecolor=red,urlcolor=Emerald,bookmarksopen,bookmarksdepth=2]{hyperref}
\usepackage{bookmark}

\theoremstyle{plain}
\newtheorem{theorem}{Theorem}[section]

\newtheorem{lemma}[theorem]{Lemma}
\newtheorem{corollary}[theorem]{Corollary}
\newtheorem{proposition}[theorem]{Proposition}
\numberwithin{equation}{section}
\theoremstyle{definition}
\newtheorem{definition}[theorem]{Definition}
\newtheorem{example}[theorem]{Example}

\providecommand{\add}{\operatorname{add}}

\providecommand{\per}{\operatorname{per}}
\providecommand{\thick}{\operatorname{thick}}

\renewcommand{\mod}{\operatorname{mod}}
\newcommand{\Hom}{\operatorname{Hom}}
\newcommand{\RHom}{\mathrm R\!\Hom}
\newcommand{\real}{\operatorname{real}}
\newcommand{\Lotimes}{\overset{\mathrm L}{\otimes}}
\newcommand{\End}{\operatorname{End}}
\newcommand{\Ext}{\operatorname{Ext}}

\newcommand{\Proj}{\operatorname{proj}}
\newcommand{\Inj}{\operatorname{inj}}

\renewcommand{\k}{\mathbf{k}}

\DeclareMathOperator{\huaEnd}{\mathcal{E}\!nd}
\DeclareMathOperator{\huaHom}{\mathcal{H}\!om}
\newcommand{\emod}[2]{#1\text{-}\!\mod #2}
\newcommand{\eheart}[2]{#2^{[-#1+1,0]}}
\newcommand{\eend}[2]{#1\text{-}\huaEnd_{#2}}
\newcommand{\PP}{\mathbf{P}}
\newcommand{\QQ}{\mathbf{Q}}
\newcommand{\RR}{\mathbf{R}}
\newcommand{\D}{\mathcal{D}}
\newcommand{\T}{\mathcal{T}}
\newcommand{\F}{\mathcal{F}}
\newcommand{\X}{\mathcal{X}}
\newcommand{\Y}{\mathcal{Y}}

\begin{document}

\title{A silting theorem for extended module categories}
\author{Chenhui Liu}
\address{LC: School of Mathematical Sciences,
Beijing Normal University,
100875 Beijing,
China}
\email{202631130047@mail.bnu.edu.cn}

\author{Yu Zhou}
\address{ZY: School of Mathematical Sciences,
Beijing Normal University,
100875 Beijing,
China}
\email{yuzhoumath@gmail.com}

\thanks{The work was supported by the National Natural Science Foundation of China (Grant No. 12671045).}

\begin{abstract}
We generalize the silting theorem of Buan and Zhou to extended module categories. For a $(d+1)$-term silting object $\mathbf{P}$ over a proper $d$-truncated non-positive dg algebra $A$, we construct a $(d+1)$-term silting object $\mathbf{Q}$ over the $d$-truncated dg endomorphism algebra $B$ of $\mathbf{P}$. The $s$-torsion pairs induced by $\mathbf{P}$ and $\mathbf{Q}$ yield equivalences of Brenner-Butler type between the torsion class on either side and the torsion-free class on the other. There is also a dg algebra homomorphism from $A$ to the $d$-truncated dg endomorphism algebra of $\mathbf{Q}$, which is a quasi-isomorphism if and only if $\operatorname{Hom}(\mathbf{P},\mathbf{P}[i])=0$ for all $i\leq-d$. Moreover, we describe Auslander-Reiten theory in ${d\text{-}\operatorname{mod}\!B}$ in terms of that in ${d\text{-}\operatorname{mod}\!A}$.
\end{abstract}

\maketitle

%==================================================%
\section{Introduction}\label{sec:introduction}
%==================================================%

Tilting theory plays a central role in the representation theory of algebras by relating module categories through torsion pairs \cite{BGP,APR,BB,HR}. It has developed into a general method for constructing and studying equivalences of categories. The derived Morita theory associated with tilting modules, complexes and bundles has important applications in the representation theory of finite groups and algebraic groups, algebraic geometry and topology; cf. \cite{Ric,AHK}. The development of cluster tilting theory has further provided categorical models for cluster algebras and their mutation combinatorics~\cite{BMRRT}. Silting complexes were introduced in~\cite{KV} as a generalization of tilting complexes. Derived Morita theory~\cite{Ke1} also applies to silting complexes, yielding equivalences with derived categories of differential graded endomorphism algebras.

The relation between 2-term silting complexes and torsion pairs was studied in~\cite{HKM,AIR}. Building on this work, the classical tilting theorem was generalized to 2-term silting complexes in~\cite{BZ}. Let $A$ be a finite-dimensional algebra over a field $\k$, and let $\mod A$ denote the category of finitely generated right $A$-modules. For a 2-term silting complex $\PP$ over $A$, let $B=\End_{D^b(\mod A)}(\PP)$. The theorem gives torsion pairs $(\T(\PP),\F(\PP))$ in $\mod A$ and $(\X(\PP),\Y(\PP))$ in $\mod B$, together with equivalences
\[\T(\PP)\simeq\Y(\PP)\text{ and }\F(\PP)\simeq\X(\PP).\]
A 2-term silting complex $\QQ$ over $B$ is also constructed, which induces the torsion pair $(\X(\PP),\Y(\PP))$. Moreover, there exists an algebra homomorphism $A\to\End_{D^b(\mod B)}(\QQ)$ which is an isomorphism if and only if $\PP$ is tilting. One key point in this theory is to realize $\mod B$ as the heart of the bounded $t$-structure on $D^b(\mod A)$ induced by $\PP$. This realization works for arbitrary silting complexes~\cite{KY}. To generalize the silting theorem, one also needs counterparts of the two torsion pairs and the second silting complex. A dg approach to the silting theorem was developed in~\cite{XYZ}. This approach also yields a partial generalization to higher-term silting objects, giving equivalences between subcategories defined by cohomological vanishing; see \cite[Section~4.5]{XYZ}. The question of whether the corresponding second object is silting after passing to the ordinary endomorphism algebra was left open.

Extended module categories provide a natural setting for a higher-term silting theorem. The correspondence between 2-term silting complexes and functorially finite torsion pairs generalizes to $(d+1)$-term silting objects and functorially finite $s$-torsion pairs in $d$-extended module categories~\cite{G,Z,WZ,HZ}. These categories also admit versions of Happel-Reiten-Smal\o{} tilting~\cite{AET,Z} and Auslander-Reiten theory~\cite{Z,MP,Ru}. Related developments concern semibricks and wide subcategories~\cite{GZ}, Auslander correspondence~\cite{M2}, higher cluster tilting~\cite{M3}, and the structure of extended hearts~\cite{M1,P,MNO}.

In this paper, we extend the silting theorem of Buan and Zhou to $(d+1)$-term silting objects by using truncated dg endomorphism algebras and extended module categories. More precisely, let $d\geq1$, and let $A$ be a proper $d$-truncated non-positive differential graded (or dg for short) $\k$-algebra, i.e., a non-positive dg algebra whose total cohomology is finite-dimensional and concentrated in $[-d+1,0]$. Let $D_{\mathrm{fd}}(A)$ be the full subcategory of $D(A)$ consisting of objects with finite-dimensional total cohomology. The $d$-extended module category of $A$ is
\[
 \emod{d}{A}
 :=\{X\in D_{\mathrm{fd}}(A)\mid
       H^i(X)=0\text{ for all }i\notin[-d+1,0]\}.
\]
Let $\PP$ be a $(d+1)$-term silting object in the perfect category $\per A$ of $A$, i.e., $\PP\in\add A*\add A[1]*\cdots*\add A[d]$, $\Hom_{D(A)}(\PP,\PP[i])=0$ for all $i>0$, and $\thick(\PP)=\per A$. The full subcategories
\[\T(\PP)=\{M\in\emod{d}{ A}\mid\Hom_{D(A)}(\PP,M[i])=0\text{ for all }i>0\}\]
and
\[\F(\PP)=\{M\in\emod{d}{ A}\mid\Hom_{D(A)}(\PP,M[i])=0\text{ for all }i\leq0\}\]
form an $s$-torsion pair in $\emod{d}{A}$ \cite[Corollary~4.1]{HZ}; see also \cite{G,Z} for the case that $A$ is an ordinary algebra. Moreover, $(\F(\PP)[d],\T(\PP))$ is an $s$-torsion pair in the $d$-extended heart $\eheart{d}{D_{\mathrm{fd}}(\PP)}$ of the bounded $t$-structure induced by $\PP$.

The $d$-truncated dg endomorphism algebra of $\PP$ is defined to be
\[B=\eend{d}{A}(\PP)=\sigma^{\geq -d+1}\sigma^{\leq 0}\left(\huaEnd_A(\PP)\right),\]
where $\huaEnd_A(\PP)$ is the dg endomorphism algebra of a $K$-projective dg $A$-module representing $\PP$, and the truncations $\sigma^{\geq -d+1}$ and $\sigma^{\leq 0}$ are canonical. For $d=1$, $B$ is isomorphic to the ordinary endomorphism algebra of $\PP$. 

The following theorem summarizes the main results (Theorem~\ref{thm:real-d-heart}, Corollary~\ref{cor:higherBB}, Proposition~\ref{prop:Qsilt}, and Theorems~\ref{thm:higherBB}, \ref{thm:equiv} and~\ref{thm:torsion-res}) of the paper.

\begin{theorem}
Let $A$ be a proper $d$-truncated non-positive dg algebra, and let $\PP$ be a $(d+1)$-term silting object in $\per A$. Let $B=\eend{d}{A}(\PP)$.
\begin{enumerate}
\item There exists an equivalence 
\[\Theta_{\PP}\colon\eheart{d}{D_{\mathrm{fd}}(\PP)}\to\emod{d}{B},\] 
which preserves and reflects triangles whose first three terms belong to the respective categories. Moreover, for any objects $M,N\in\eheart{d}{D_{\mathrm{fd}}(\PP)}$, there exist natural isomorphisms
\[\Hom_{D(A)}\left(M,N[i]\right)\cong\Hom_{D(B)}(\Theta_{\PP}(M),\Theta_{\PP}(N)[i]),\quad i\leq1.\]
\item Let $\X(\PP)=\Theta_{\PP}(\F(\PP)[d])$ and $\Y(\PP)=\Theta_{\PP}(\T(\PP))$. Then $(\X(\PP),\Y(\PP))$ is an $s$-torsion pair in $\emod{d}{B}$, with equivalences 
\[\T(\PP)\simeq\Y(\PP)\text{ and }\F(\PP)\simeq\X(\PP).\]
\item Let $\Gamma=\sigma^{\leq 0}\huaEnd_A(\PP)$ be the non-positive truncation of $\huaEnd_A(\PP)$. Then
\[\QQ:=\RHom_A(\PP,A[d])\Lotimes_{\Gamma}B\] 
is a $(d+1)$-term silting object in $\per B$, whose induced $s$-torsion pair $(\T(\QQ),\F(\QQ))$ satisfies $\T(\QQ)=\X(\PP)$ and $\F(\QQ)=\Y(\PP)$.

\item Let $C=\eend{d}{B}(\QQ)$, using a $K$-projective dg $A$-$B$-bimodule representing $\QQ$. There exists a dg algebra homomorphism 
\[\vartheta\colon A\to C,\] 
which is a quasi-isomorphism if and only if $\Hom_{\per A}(\PP,\PP[i])=0$ for all $i\leq-d$.
\item Let $\vartheta_*\colon\emod{d}{C}\to\emod{d}{A}$ be the restriction functor induced by $\vartheta$. Then  $\vartheta_*$ restricts to equivalences $\X(\QQ)\xrightarrow{\simeq}\T(\PP)$ and $\Y(\QQ)\xrightarrow{\simeq}\F(\PP)$.
\end{enumerate}
\end{theorem}

When $d=1$ and $A$ is an ordinary finite-dimensional algebra, this theorem recovers the silting theorem of Buan and Zhou. In that case, the functor $\vartheta_*$ is fully faithful. For $d>1$, this may fail; see Example~\ref{ex:A3-worked}. 

If $(\T(\PP),\F(\PP))$ is split, then $\vartheta$ is a quasi-isomorphism. More precisely, Proposition~\ref{prop:separating} bounds the negative self-extensions of $\PP$ in terms of the cohomological range of $A$. In particular, if $A$ is an ordinary finite-dimensional algebra and $(\T(\PP),\F(\PP))$ is split, then $\PP$ is tilting, generalizing \cite[Proposition~5.7]{BZ}.

We apply the higher Brenner-Butler theorem to Auslander-Reiten theory in $\emod{d}{B}$. We describe the connecting almost split triangles and, when $(\X(\PP),\Y(\PP))$ is split, determine all almost split triangles in $\emod{d}{B}$; see Propositions~\ref{prop:connecting} and \ref{prop:AR-splitting}.

The paper is organized as follows. Section~\ref{sec:pre} recalls preliminaries on dg algebras, derived categories, adjunctions and extended module categories. Section~\ref{sec:s-tor} discusses $s$-torsion pairs and almost split triangles. In Section~\ref{sec:higherBB}, we prove the higher Brenner-Butler theorem and construct the silting object $\QQ$. Section~\ref{sec:algC} studies the $d$-truncated dg endomorphism algebra $C$ of $\QQ$ and the dg algebra morphism $\vartheta\colon A\to C$. Section~\ref{sec:AR-theory} applies the theorem to Auslander-Reiten theory. Section~\ref{sec:sep-spli} is devoted to separating and splitting silting objects.

%==================================================%
\subsection*{Convention}
%==================================================%

Throughout this paper, $\k$ is a field, and $d$ is a positive integer. We write $\mathbb D=\Hom_{\k}(-,\k)$ for the $\k$-linear dual. Any subcategory of a category is assumed to be full and closed under isomorphisms. We use $X\in\X$ to denote that $X$ is an object in a category $\X$, and use $\X\subseteq\Y$ to denote that $\X$ is a subcategory of $\Y$. For an object $X$ of an additive category, we denote by $\add X$ the full subcategory consisting of direct summands of finite direct sums of copies of $X$. For any two morphisms $f\colon Y\to Z$ and $g\colon X\to Y$, we denote by $f\circ g\colon X\to Z$ their composition.

Let $\Y$ be a category. For any $\X\subseteq\Y$ and any $Y\in\Y$, a right $\X$-approximation of $Y$ is a morphism $f:X\to Y$ with $X\in\X$ and such that any morphism from an object $X'$ in $\X$ to $Y$ factors through $f$. A left $\X$-approximation of $Y$ is defined dually. A subcategory $\X$ of $\Y$ is called contravariantly (resp. covariantly) finite in $\Y$ if any object in $\Y$ has a right (resp. left) $\X$-approximation. A subcategory $\X$ of $\Y$ is called functorially finite if it is both contravariantly finite and covariantly finite. A morphism $f:X\to Y$ is called right minimal if every morphism $g:X\to X$ such that $f\circ g=f$ is an isomorphism. A left minimal morphism is defined dually.

The shift functor in a triangulated category is denoted by $[1]$. For any two subcategories $\X$ and $\Y$ of a triangulated category, we denote by $\X*\Y$ the full subcategory consisting of objects $Z$ such that there is a triangle
\[X\to Z\to Y\to X[1],\]
with $X\in\X$ and $Y\in\Y$. We use $\Hom(\X,\Y)=0$ to denote $\Hom(X,Y)=0$ for any $X\in\X$ and $Y\in\Y$. A subcategory $\X$ of a triangulated category is called closed under extensions if $\X*\X\subseteq\X$. 

Let $\X$ be an extension-closed subcategory of a triangulated category $\D$. A triangle in $\X$ means a triangle in $\D$ whose first three terms belong to $\X$. An object $X$ in $\X$ is called $\Ext$-projective (resp. $\Ext$-injective) if $\Hom_{\D}(X,Y[1])=0$ (resp. $\Hom_{\D}(Y,X[1])=0$) for all $Y\in\X$. We denote by $\Proj\X$ and $\Inj\X$ the subcategories of $\X$ consisting of $\Ext$-projective and $\Ext$-injective objects in $\X$, respectively. We say that $\X$ has enough projectives (resp. injectives) if for any object $X\in\X$, there exists a triangle in $\X$:
\[Y\to P\to X\to Y[1]\text{ (resp. $X\to I\to Z\to X[1]$)},\]
with $P\in\Proj\X$ (resp. $I\in\Inj\X$).

%==================================================%
\section{Preliminaries}\label{sec:pre}
%==================================================%

%==================================================%
\subsection{Differential graded algebras, derived categories and adjunctions}\label{subsec:dg-pre}\ 
%==================================================%

We recall some basic notions and results on differential graded (dg) algebras, derived categories and adjunctions; cf. \cite{Ke1,Ke2,Kr,Rie,Y}.

All dg algebras and dg modules are defined over $\k$ and are $\mathbb Z$-graded, with differentials of degree $1$. Their underlying dg $\k$-modules are complexes of $\k$-vector spaces. For $n\in\mathbb Z$, we denote by $\sigma^{\leq n}$ and $\sigma^{\geq n}$ the canonical truncations of complexes of $\k$-vector spaces; cf. \cite[Appendix~A]{Z}. We write $H^i(X)$ for the $i$-th cohomology of a complex of $\k$-vector spaces $X$. Thus, we have
\[\sigma^{\geq i}(\sigma^{\leq i}(X))\cong H^i(X)[-i]\cong \sigma^{\leq i}(\sigma^{\geq i}(X)).\]

Unless otherwise stated, all dg modules are right dg modules. We write ${}_RM_S$ for a dg bimodule with the left $R$-action and right $S$-action. We omit action subscripts already determined by the ambient category or the subscript of $\Hom$ or $\otimes$.

Let $\Lambda$ be a dg algebra. For dg $\Lambda$-modules $X,Y$, we denote the dg Hom complex by $\huaHom_\Lambda(X,Y)$ and the dg endomorphism algebra by $\huaEnd_\Lambda(X)=\huaHom_\Lambda(X,X)$, with multiplication given by composition. We denote by $D(\Lambda)$ the derived category of right dg $\Lambda$-modules, by $\per \Lambda$ its perfect derived category, and by $D_{\mathrm{fd}}(\Lambda)$ the full subcategory of objects with finite-dimensional total cohomology.

A dg $\Lambda$-module $P$ is \emph{$K$-projective} if $\huaHom_\Lambda(P,N)$ is acyclic for every acyclic dg $\Lambda$-module $N$. A \emph{$K$-projective resolution} of a dg $\Lambda$-module $X$ is a quasi-isomorphism $P_X\to X$ where $P_X$ is $K$-projective. Every dg module admits such a resolution. For dg $\Lambda$-modules $X,N$, the derived Hom $\RHom_\Lambda(X,N)\in D(\k)$ is represented by $\huaHom_\Lambda(P_X,N)$, where $P_X\to X$ is a $K$-projective resolution. For every $n\in\mathbb Z$, we have
\[
 H^n\RHom_\Lambda(X,N)
 \cong\Hom_{D(\Lambda)}(X,N[n]).
\]
Let $\Lambda'$ be another dg algebra and let ${}_\Lambda T_{\Lambda'}$ be a dg $\Lambda$-$\Lambda'$-bimodule. For a dg $\Lambda$-module $X$, the derived tensor product $X\Lotimes_\Lambda T\in D(\Lambda')$ is represented by $P_X\otimes_\Lambda T$. If, moreover, $X$ is a dg $R$-$\Lambda$-bimodule for some dg algebra $R$, we choose a $K$-projective resolution $P_X\to X$ in the category of dg $R^{\mathrm{op}}\otimes_{\k}\Lambda$-modules. Since $P_X$ remains $K$-projective as a right $\Lambda$-module, the above constructions then give $\RHom_\Lambda(X,N)\in D(R)$ and $X\Lotimes_\Lambda T\in D(R^{\mathrm{op}}\otimes_{\k}\Lambda')$. More generally, if $N$ is a dg $S$-$\Lambda$-bimodule for some dg algebra $S$, then $\RHom_\Lambda(X,N)$ belongs to $D(S^{\mathrm{op}}\otimes_{\k}R)$.

For $M\in D(\Lambda)$ and $N\in D(\Lambda')$, there is a natural isomorphism in $D(\k)$:
\begin{equation}\label{eq:tensor-Hom}
 \RHom_{\Lambda'}
   \left(M\overset{\mathrm L}{\otimes}_{\Lambda}T,N\right)
 \cong
 \RHom_{\Lambda}
   \left(M,\RHom_{\Lambda'}({}_\Lambda T,N)\right).
\end{equation}
Taking $H^0$ provides the adjunction $-\overset{\mathrm L}{\otimes}_\Lambda T\dashv\mathrm R\!\Hom_{\Lambda'}({}_\Lambda T,-)$ between $D(\Lambda)$ and $D(\Lambda')$. Its unit and counit are
\begin{equation}\label{eq:unit}
    \eta_M\colon M\rightarrow\RHom_{\Lambda'}\left({}_\Lambda T,M\overset{\mathrm L}{\otimes}_\Lambda T\right),
\end{equation}
\begin{equation}\label{eq:counit}
    \varepsilon_N\colon \RHom_{\Lambda'}({}_\Lambda T,N)\overset{\mathrm L}{\otimes}_\Lambda T\rightarrow N
\end{equation}
in $D(\Lambda)$ and $D(\Lambda')$, respectively.

More generally, if ${}_RM_\Lambda$ and ${}_SN_{\Lambda'}$ are dg bimodules, then \eqref{eq:tensor-Hom} is an isomorphism in $D(S^{\mathrm{op}}\otimes_{\k}R)$. The unit and counit in \eqref{eq:unit} and \eqref{eq:counit} are morphisms in $D(R^{\mathrm{op}}\otimes_{\k}\Lambda)$ and $D(S^{\mathrm{op}}\otimes_{\k}\Lambda')$, respectively. Whenever this adjunction is an equivalence, the morphisms in \eqref{eq:unit} and \eqref{eq:counit} are isomorphisms.

For a dg algebra homomorphism $f\colon\Lambda\to\Lambda'$, the right $\Lambda'$-module $\Lambda'$ becomes a dg $\Lambda$-$\Lambda'$-bimodule ${}_\Lambda\Lambda'_{\Lambda'}$, with left $\Lambda$-action induced by $f$. We denote by 
\[f_*\colon D(\Lambda')\to D(\Lambda)\]
the restriction functor along $f$, and denote by
\[f^*\colon D(\Lambda)\to D(\Lambda')\]
the functor given by $f^*(M)=M\Lotimes_{\Lambda}\Lambda'$. The morphism $f$ is a quasi-isomorphism if and only if $f_*$ is an equivalence.

For any $N\in D(\Lambda')$, there exists a natural isomorphism in $D(\Lambda)$:
\begin{equation}\label{eq:HomN=N}
    \mathrm R\!\Hom_{\Lambda'}({}_\Lambda\Lambda',N)
\cong f_*(N).
\end{equation}
Hence, we obtain the adjunction
\[
f^*\dashv f_*.
\]
In this case, using \eqref{eq:HomN=N}, the unit and counit in \eqref{eq:unit} and \eqref{eq:counit} become
\begin{equation}\label{eq:tensor-unit}
 \eta_M\colon M\cong M\overset{\mathrm L}{\otimes}_\Lambda\Lambda\xrightarrow{1_M\overset{\mathrm L}{\otimes}_\Lambda f} f_*f^*(M),
\end{equation}
\begin{equation}\label{eq:tensor-counit}
 \varepsilon_N\colon f^*f_*(N)\to N,
\end{equation}
and the natural isomorphism given by this adjunction is
\begin{equation}\label{eq:adjun-exp}
  \begin{array}{ccc}
    \Hom_{D(\Lambda')}(f^*(M),N) & \xrightarrow{\cong} & \Hom_{D(\Lambda)}(M,f_*(N))\\
    h & \mapsto & f_*(h)\circ\eta_M
  \end{array}.
\end{equation}
For the case $M=\Lambda$, under $f^*(\Lambda)\cong\Lambda'$, we have
\begin{equation}\label{eq:unit-algebra}
 \eta_\Lambda=f\colon\Lambda\to f_*(\Lambda').
\end{equation}

A dg algebra $\Lambda$ is \emph{proper} if its total cohomology is finite dimensional. It is \emph{non-positive} if $\Lambda^i=0$ for $i>0$, and \emph{connective} if $H^i(\Lambda)=0$ for $i>0$. Following \cite[Definition~1.1]{M2}, it is \emph{$d$-truncated} if $H^i(\Lambda)=0$ for $i\leq-d$. Thus, a connective $d$-truncated dg algebra has cohomology concentrated in $[-d+1,0]$. 

If $\Lambda$ is non-positive, then $D(\Lambda)$ has the canonical $t$-structure $(D(\Lambda)^{\leq0},D(\Lambda)^{\geq0})$, where, for $n\in\mathbb Z$,
\[D(\Lambda)^{\leq n}=\{X\in D(\Lambda)\mid H^i(X)=0\text{ for all }i>n\},\]
\[
D(\Lambda)^{\geq n}=\{X\in D(\Lambda)\mid H^i(X)=0\text{ for all }i<n\}.
\]
This $t$-structure restricts to a bounded $t$-structure on $D_{\mathrm{fd}}(\Lambda)$. The canonical truncations above induce the truncation functors $\sigma_\Lambda^{\leq n}$ and $\sigma_\Lambda^{\geq n}$ on $D(\Lambda)$. For each $X\in D(\Lambda)$, they give a triangle in $D(\Lambda)$
\[
 \sigma_\Lambda^{\leq n}(X)\to X\to\sigma_\Lambda^{\geq n+1}(X)\to\sigma_\Lambda^{\leq n}(X)[1].
\]
For a morphism $f\colon\Lambda\to\Lambda'$ of non-positive dg algebras, restriction preserves cohomology, is standard $t$-exact, commutes with the canonical truncations, and reflects isomorphisms. 

%==================================================%
\subsection{Extended hearts and extended module categories}\ 
%==================================================%

Let $(\D^{\leq0},\D^{\geq0})$ be a bounded $t$-structure on a triangulated category $\D$. Its \emph{$d$-extended heart} is
\[
 \eheart{d}{\D}
 :=\D^{\geq-d+1}\cap\D^{\leq0}.
\]
It is closed under extensions, and hence an extriangulated category in the sense of \cite{NP}. When $d=1$, the $d$-extended heart coincides with the classical heart.

For a non-positive dg algebra $\Lambda$, its \emph{$d$-extended module category} is the $d$-extended heart of the standard $t$-structure on $D_{\mathrm{fd}}(\Lambda)$:
\[
 \emod{d}{\Lambda}:=D_{\mathrm{fd}}(\Lambda)^{[-d+1,0]}=\{X\in D_{\mathrm{fd}}(\Lambda)\mid H^i(X)=0\ \forall i\notin [-d+1,0]\}.
\]
If $\Lambda$ is proper, this category is Hom-finite, idempotent complete and Krull-Schmidt (see \cite{MP}). 

We compare the $d$-extended module categories of a non-positive dg algebra and its $d$-truncated version. This is also mentioned in \cite[Remark~1.10]{M2}.

\begin{proposition}\label{prop:trun-equiv} 
  Let $\Lambda$ be a proper non-positive dg algebra. Let $\Lambda'=\sigma^{\geq -d+1}(\Lambda)$ and let $f\colon \Lambda\to\Lambda'$ be the canonical projection. Then the restriction functor $f_*\colon D(\Lambda')\to D(\Lambda)$ induces an equivalence
  \[
  f_*\colon\emod{d}{\Lambda'}\xrightarrow{\simeq}\emod{d}{\Lambda},
  \]
  with quasi-inverse $X\mapsto\sigma_{\Lambda'}^{\leq0}\RHom_\Lambda({}_{\Lambda'}\Lambda'_\Lambda,X)$, where $\Lambda'$ is regarded as a right dg $\Lambda$-module via $f$. Moreover, for all $X,Y\in\emod{d}{\Lambda'}$, we have natural isomorphisms
  \begin{equation}\label{eq:isoleq1}
    \Hom_{D(\Lambda')}(X,Y[i])\cong\Hom_{D(\Lambda)}(f_*X,f_*Y[i]),\quad i\leq 1.
  \end{equation}
\end{proposition}

\begin{proof}
Let $f^!=\RHom_\Lambda({}_{\Lambda'}\Lambda'_\Lambda,-)\colon D(\Lambda)\to D(\Lambda')$ be the right adjoint of $f_*$. For any $X\in\emod{d}{\Lambda}$, we have $f^!(X)\in D(\Lambda')^{\geq-d+1}$. Moreover, the vector space
\[H^j(f^!(X))\cong\Hom_{D(\Lambda')}(\Lambda',f^!(X)[j])\cong\Hom_{D(\Lambda)}(f_*\Lambda',X[j])\]
is finite-dimensional. Hence, $\sigma_{\Lambda'}^{\leq 0}f^!(X)\in\emod{d}{\Lambda'}$. Applying $\mathrm R\!\Hom_\Lambda(-,X)$ to the exact sequence of dg $\Lambda$-$\Lambda$-bimodules $0\to \ker f\to \Lambda\xrightarrow{f}\Lambda'\to 0$ gives a triangle in $D(\Lambda)$:
\begin{equation}\label{eq:bi-tri}
    f_*f^!(X)\xrightarrow{\alpha_X}X\to\RHom_\Lambda(\ker f,X)\to\left(f_*f^!(X)\right)[1],
\end{equation}
where $\alpha$ is the counit of $f_*\dashv f^!$. Since $\ker f\in D(\Lambda)^{\leq-d}$ and $X\in D(\Lambda)^{\geq-d+1}$, we have $H^j\mathrm R\!\Hom_\Lambda(\ker f,X)\cong\Hom_{D(\Lambda)}(\ker f,X[j])=0$ for all $j\leq0$. Consequently, $\alpha_X$ induces isomorphisms in all degrees $j\leq0$. It follows that there are natural isomorphisms
\[
 f_*\!\left(\sigma_{\Lambda'}^{\leq0}f^!(X)\right) \cong \sigma_{\Lambda}^{\leq0}\!\left(f_*f^!(X)\right)\xrightarrow[\cong]{\sigma^{\leq 0}_{\Lambda}(\alpha_X)} X,
\]
which implies the natural isomorphism between $f_*\circ\sigma_{\Lambda'}^{\leq0}\circ f^!$ and the identity functor on $\emod{d}{\Lambda}$. In particular, $f_*$ is dense.

Since $H^1(X)=0$, the triangle~\eqref{eq:bi-tri} also implies $H^1(f_*f^!(X))=0$ and hence $H^1(\alpha_X)$ is also an isomorphism.

For $V\in\emod{d}{\Lambda'}$, let $\beta_V\colon V\to f^!f_*(V)$ be the unit. The adjunction identity
\[
 \alpha_{f_*(V)}\circ f_*(\beta_V)=1_{f_*(V)}
\]
and the isomorphism $\sigma_{\Lambda}^{\leq1}(\alpha_{f_*(V)})$ imply that $\sigma_{\Lambda'}^{\leq1}(\beta_V)$ is an isomorphism. Consequently, we have $\operatorname{cone}(\beta_V)\in D(\Lambda')^{\geq 2}$. Therefore, for any $U\in\emod{d}{\Lambda'}$ and $i\leq1$, we have natural isomorphisms
\[
 \Hom_{D(\Lambda')}(U,V[i])\xrightarrow[\cong]{\,\beta_V[i]\circ(-)\,}\Hom_{D(\Lambda')}(U,f^!f_*(V)[i])\cong\Hom_{D(\Lambda)}(f_*(U),f_*(V)[i]).
\]
The composition of these two isomorphisms sends $g\colon U\to V[i]$ to $\alpha_{f_*(V)}[i]\circ f_*\left(\beta_V[i]\circ g\right)=f_*(g)$. Hence, we obtain the required natural isomorphisms. Taking $i=0$ shows that $f_*$ is fully faithful on $\emod{d}{\Lambda'}$. Hence, $f_*$ is an equivalence, with quasi-inverse $\sigma_{\Lambda'}^{\leq0}\circ f^!$.
\end{proof}

We next consider restriction along morphisms between $d$-truncated dg algebras.

\begin{proposition}\label{prop:trun}
Let $\Lambda$ and $\Lambda'$ be proper $d$-truncated non-positive dg algebras, and let $f\colon\Lambda\to\Lambda'$ be a dg algebra homomorphism. Then $f$ is a quasi-isomorphism if and only if $f_*$ restricts to an equivalence:
\[f_*\colon\emod{d}{\Lambda'}\to\emod{d}{\Lambda}.\]
\end{proposition}

\begin{proof}
Suppose that $f$ is a quasi-isomorphism. Then $f_*\colon D(\Lambda')\to D(\Lambda)$ is an equivalence. Since $f_*$ preserves cohomology, it restricts to the required equivalence.

Conversely, suppose that $f_*\colon\emod{d}{\Lambda'}\to\emod{d}{\Lambda}$ is an equivalence. By \eqref{eq:unit-algebra}, we have $\eta_\Lambda=f\colon\Lambda\to f_*(\Lambda')$ under $f^*(\Lambda)\cong\Lambda'$. For any object $Y\in\emod{d}{\Lambda'}$, by \eqref{eq:adjun-exp}, there exists a natural isomorphism
\[
  \begin{array}{ccc}
    \Hom_{D(\Lambda')}(\Lambda',Y) & \xrightarrow{\cong} & \Hom_{D(\Lambda)}(\Lambda,f_*(Y))\\
    h & \mapsto & f_*(h)\circ f
  \end{array}.
\]
Since $f_*$ is fully faithful, this induces a natural isomorphism
\[
  \begin{array}{ccc}
    \Hom_{D(\Lambda)}(f_*(\Lambda'),f_*(Y)) & \xrightarrow{\cong} & \Hom_{D(\Lambda)}(\Lambda,f_*(Y))\\
    g & \mapsto & g\circ f
  \end{array}.
\]
Since $f_*$ is dense, we can replace $f_*(Y)$ in the above natural isomorphism by any object in $\emod{d}{\Lambda}$. Since $\Lambda$ and $\Lambda'$ are proper and $d$-truncated, we have $\Lambda\in\emod{d}{\Lambda}$ and $\Lambda'\in\emod{d}{\Lambda'}$. Thus, by the Yoneda Lemma, $f\colon\Lambda\to f_*(\Lambda')$ is an isomorphism in $\emod{d}{\Lambda}$. Hence $f$ is a quasi-isomorphism.
\end{proof}

For the remainder of this subsection, let $\Lambda$ be a proper $d$-truncated non-positive dg algebra. By \cite[Corollary~2.3]{MP}, $\emod{d}{\Lambda}$ has enough projectives and injectives, with
\begin{equation}\label{eq:extend-proj-inj}
 \Proj\emod{d}{\Lambda}=\add \Lambda \text{ and } \Inj\emod{d}{\Lambda}=\add\left(\mathbb D(\Lambda)[d-1]\right).
\end{equation}

By taking cohomology, we may relate the extended module category $\emod{d}{\Lambda}$ with the category $\mod H^0(\Lambda)$ of finitely generated right $H^0(\Lambda)$-modules as follows.
\begin{enumerate}
  \item There exist equivalences
    \begin{equation}\label{eq:proj-equiv}
        H^0:\add\Lambda\xrightarrow{\simeq}\operatorname{proj}H^0(\Lambda),
    \end{equation}
    \begin{equation}\label{eq:inj-equiv}
        H^{-d+1}:\add(\mathbb D\Lambda[d-1])\xrightarrow{\simeq}\operatorname{inj}H^0(\Lambda),
    \end{equation}
    where $\operatorname{proj}H^0(\Lambda)$ (resp., $\operatorname{inj}H^0(\Lambda)$) denotes the subcategory of $\mod H^0(\Lambda)$ consisting of projective (resp. injective) objects.
  \item For any morphism $u:X\to Y$ in $\emod{d}{\Lambda}$, we have equivalences
  \begin{equation}\label{eq:cone-lower-bound}
    \operatorname{cone}(u)\in\emod{d}{\Lambda} \iff H^{-d+1}(u)\text{ is injective},  
  \end{equation}
  \begin{equation}\label{eq:cone-upper-bound}
      \operatorname{cone}(u)[-1]\in\emod{d}{\Lambda} \iff H^0(u)\text{ is surjective}.
  \end{equation}
Cf. also \cite[Remark~1.1]{MP}.
  \item For any $X\in\emod{d}{\Lambda}$ and any $S\in\mod H^0(\Lambda)$, regarding $S$ as a dg $\Lambda$-module via $\Lambda\to H^0(\Lambda)$, we have natural isomorphisms
  \begin{equation}\label{eq:simple-top-Hom}
      \Hom_{D(\Lambda)}(X,S)\cong\Hom_{H^0(\Lambda)}(H^0(X),S),
  \end{equation}
  \begin{equation}\label{eq:simple-socle-Hom}
      \Hom_{D(\Lambda)}(S[d-1],X)\cong\Hom_{H^0(\Lambda)}(S,H^{-d+1}(X)).
  \end{equation}
\end{enumerate}

\begin{definition}[{\cite[Definition~3.6]{MP}}]\label{def:rad-soc}
For any object $X\in\emod{d}{\Lambda}$, we define
\[
 \operatorname{top}X=\operatorname{top}_{H^0(\Lambda)}H^0(X) \text{ and } \operatorname{soc}X=\left(\operatorname{soc}_{H^0(\Lambda)}H^{-d+1}(X)\right)[d-1],
\]
where the modules on the right are regarded as dg $\Lambda$-modules via $\Lambda\to H^0(\Lambda)$. Define $\operatorname{rad}X$ and $X/\operatorname{soc}X$ by the following triangles:
\[
\operatorname{rad}X\to X\to\operatorname{top}X\to(\operatorname{rad}X)[1],
\]
\[
\operatorname{soc}X\to X\to X/\operatorname{soc}X\to(\operatorname{soc}X)[1].
\]
Both $\operatorname{rad}X$ and $X/\operatorname{soc}X$ belong to $\emod{d}{\Lambda}$ by the equivalences \eqref{eq:cone-lower-bound} and \eqref{eq:cone-upper-bound}.
\end{definition}

%==================================================%
\section{\texorpdfstring{$s$}{s}-torsion pairs and almost split triangles}\label{sec:s-tor}
%==================================================%

In this section, we study basic properties of $s$-torsion pairs and almost split triangles in extended module categories. We begin in the more general setting of extended hearts.

Let $\eheart{d}{\D}$ be the $d$-extended heart of a bounded $t$-structure on a triangulated category $\D$. An \emph{$s$-torsion pair} (\cite[Definition~3.1]{AET}) in $\eheart{d}{\D}$ is a pair $(\T,\F)$ of full subcategories closed under isomorphisms such that 
\begin{itemize}
    \item $\Hom_{\D}(\T,\F)=0$,
    \item $\eheart{d}{\D}=\T*\F$, and
    \item $\Hom_{\D}(\T,\F[-1])=0$.
\end{itemize}
The last condition implies $\Hom_{\D}(\T,\F[i])=0$ for all $i<0$ by \cite[Corollary~1.10]{Z}. For any $s$-torsion pair $(\T,\F)$ in $\eheart{d}{\D}$, we have 
\[\T=\{X\in\eheart{d}{\D}\mid \Hom_{\D}(X,\F)=0\},\]
\[\F=\{X\in\eheart{d}{\D}\mid \Hom_{\D}(\T,X)=0\}.\]
Moreover, any object $X\in\eheart{d}{\D}$ admits a triangle
\begin{equation}\label{eq:torsion-triangle}
 \mathfrak{t}X\to X\to \mathfrak{f}X\to(\mathfrak{t}X)[1],
\end{equation}
with $\mathfrak{t}X\in\T$ and $\mathfrak{f}X\in\F$. This triangle is unique up to isomorphism (\cite[Proposition 3.7]{AET}), and is called the canonical triangle of $X$ with respect to the $s$-torsion pair $(\T,\F)$. These triangles define additive functors
\begin{equation}\label{eq:tor-functors}
    \mathfrak{t}:\eheart{d}{\D}\to\T\text{ and }\mathfrak{f}:\eheart{d}{\D}\to\F,
\end{equation}
which are called \emph{torsion} and \emph{torsion-free} functors with respect to the $s$-torsion pair $(\T,\F)$.

Let $A$ be a proper $d$-truncated non-positive dg algebra.

\begin{lemma}\label{lem:torsion-Ext-inj}
Let $(\T,\F)$ be an $s$-torsion pair in
$\emod{d}{A}$. Then
\[
 \Inj\T
 =\add\mathfrak{t}(\mathbb D A[d-1])\text{ and }
 \Proj\F=\add\mathfrak{f} A.
\]
Moreover, $\T$ has enough Ext-injectives and $\F$ has enough Ext-projectives.
\end{lemma}

\begin{proof}
  We only show the assertions for $\T$, because those for $\F$ can be shown dually.

  For any $I\in\add\left(\mathbb D A[d-1]\right)$ and $T\in\T$, applying $\Hom_{D(A)}(T,-)$ to the canonical triangle of $I$ with respect to $(\T,\F)$ yields an exact sequence
  \[
  \Hom_{D(A)}(T,\mathfrak{f}I)\to\Hom_{D(A)}(T,\mathfrak{t}I[1])\to\Hom_{D(A)}(T,I[1]),
  \]
  where the first term vanishes since $T\in\T$ and $\mathfrak{f}I\in\F$, and the last term vanishes, since by \eqref{eq:extend-proj-inj}, $I\in\Inj\emod{d}{A}$. Hence, we have $\Hom_{D(A)}(T,\mathfrak{t}I[1])=0$, which implies $\mathfrak{t}I\in\Inj\T$.

  For any $X\in\T$, since $\emod{d}{A}$ has enough injectives $\add\left(\mathbb D A[d-1]\right)$, there exists $I_0\in \add\left(\mathbb D A[d-1]\right)$ and a triangle $X\xrightarrow{\alpha}I_0\to L\to X[1]$ with $L\in\emod{d}{A}$. Consider the canonical triangle of $I_0$ with respect to $(\T,\F)$:
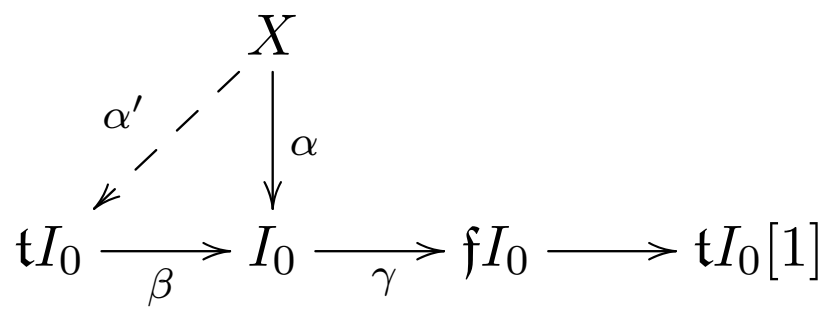
  \[\xymatrix{
    &X\ar[d]^{\alpha}\ar@{-->}[ld]_{\alpha'}\\
    \mathfrak{t}I_0\ar[r]_{\beta} &I_0\ar[r]_{\gamma} & \mathfrak{f}I_0\ar[r]& \mathfrak{t}I_0[1]
  }\]
  We have $\gamma\circ\alpha=0$ by $X\in\T$ and $\mathfrak{f}(I_0)\in\F$. So there exists a morphism $\alpha'\colon X\to \mathfrak{t}I_0$ such that $\alpha=\beta\circ\alpha'$. The injectivity of $H^{-d+1}(\alpha)$ implies that of $H^{-d+1}(\alpha')$. Hence, by \eqref{eq:cone-lower-bound}, $\operatorname{cone}(\alpha')\in\emod{d}{A}$. Since both $X$ and $\mathfrak{t}I_0$ belong to $\T$, by \cite[Lemma~3.9]{GZ}, it follows that $\operatorname{cone}(\alpha')\in\T$. Consequently, we obtain a triangle
  \[
  X\xrightarrow{\alpha'}\mathfrak{t}I_0\to\operatorname{cone}(\alpha')\to X[1]
  \]
  in $\T$ with middle term Ext-injective, which implies that $\T$ has enough injectives. If $X\in\Inj\T$, this triangle splits. Thus, we have $\Inj\T=\add\mathfrak{t}(\mathbb D A[d-1])$. 
\end{proof}

An object $\PP\in\per A$ is called \emph{silting} if $\Hom_{\per A}(\PP,\PP[i])=0$ for all $i>0$, and $\thick(\PP)=\per A$. It is called \emph{tilting} if, in addition, $\Hom_{D(A)}(\PP,\PP[i])=0$ for all $i<0$. An object of $\per A$ is called \emph{$(d+1)$-term} if it belongs to $\add A*\add A[1]*\cdots*\add A[d]$.

Let $\PP$ be a $(d+1)$-term silting object in $\per A$. We will use the following consequence of \cite[Lemma~3.7]{AMY}.
\begin{lemma}\label{lem:AinP}
We have $A[d]\in\add\PP*\add\PP[1]*\cdots*\add\PP[d]$.
\end{lemma}

By \cite[Proposition~2.1]{WZ}, the subcategories
\[
 D_{\mathrm{fd}}(\PP)^{\leq0}
 =\{X\in D_{\mathrm{fd}}(A)\mid
       \Hom_{D(A)}(\PP,X[i])=0\text{ for }i>0\},
\]
\[
 D_{\mathrm{fd}}(\PP)^{\geq0}
 =\{X\in D_{\mathrm{fd}}(A)\mid
       \Hom_{D(A)}(\PP,X[i])=0\text{ for }i<0\}
\]
form a bounded $t$-structure, and we have
\begin{equation}\label{eq:silt-inc}
\begin{aligned}
 D_{\mathrm{fd}}(A)^{\leq-d}
 &\subseteq D_{\mathrm{fd}}(\PP)^{\leq0}
 \subseteq D_{\mathrm{fd}}(A)^{\leq0},\\
 D_{\mathrm{fd}}(A)^{\geq1}
 &\subseteq D_{\mathrm{fd}}(\PP)^{\geq1}
 \subseteq D_{\mathrm{fd}}(A)^{\geq-d+1}.
\end{aligned}
\end{equation}
Then the $d$-extended heart of this $t$-structure is
\[\begin{array}{rcl}
  \eheart{d}{D_{\mathrm{fd}}(\PP)} & = & D_{\mathrm{fd}}(\PP)^{\leq 0}\cap D_{\mathrm{fd}}(\PP)^{\geq -d+1}\\
  & = & \{X\in D_{\mathrm{fd}}(A)\mid \Hom_{D(A)}(\PP,X[i])=0\text{ for }i\notin[-d+1,0]\}.
\end{array}
\]
The subcategories
\[
 \T(\PP)=D_{\mathrm{fd}}(\PP)^{\leq0}\cap\emod{d}{A}\text{ and }
 \F(\PP)=D_{\mathrm{fd}}(\PP)^{\geq1}\cap\emod{d}{A}
\]
form a functorially finite $s$-torsion pair in $\emod{d}{A}$; see \cite[Corollary~4.1]{HZ}. It is straightforward to check that 
\[\eheart{d}{D_{\mathrm{fd}}(\PP)}\cap\emod{d}{A}=\eheart{d}{D_{\mathrm{fd}}(\PP)}\cap D(A)^{\geq -d+1}=\T(\PP),\]
\[\eheart{d}{D_{\mathrm{fd}}(\PP)}\cap\emod{d}{A}[d]=\eheart{d}{D_{\mathrm{fd}}(\PP)}\cap D(A)^{\leq -d}=\F(\PP)[d].\]
We denote by $\mathfrak{t}_{\PP}$ and $\mathfrak{f}_{\PP}$ the torsion and torsion-free functors with respect to the $s$-torsion pair $(\T(\PP),\F(\PP))$; see \eqref{eq:tor-functors}.

The Nakayama functor $\nu_A=-\Lotimes_A\mathbb D A$ restricts to an equivalence $\add A\simeq\add\mathbb D A$. For $P\in\per A$ and $M\in D_{\mathrm{fd}}(A)$, there are natural isomorphisms
\begin{equation}\label{eq:Serre}
 \mathbb D\Hom_{D(A)}(P,M)
 \cong\Hom_{D(A)}(M,\nu_A P).
\end{equation}

The following proposition gives the dg version of \cite[Proposition~2.14 and Corollary~4.13]{Z}. For an object $X$, let $|X|$ denote the number of isomorphism classes of its indecomposable direct summands.

\begin{proposition}\label{prop:Ext-proj-inj}
The following hold.
\begin{enumerate}
\item Both $\T(\PP)$ and $\F(\PP)$ have enough Ext-projectives and Ext-injectives, and
\begin{align*}
    \Proj\T(\PP)&=\add \sigma_A^{\geq-d+1}(\PP),&\Inj\T(\PP)&=\add\mathfrak{t}_{\PP}(\mathbb D A[d-1]),\\
    \Proj\F(\PP)&=\add\mathfrak{f}_{\PP}(A),&\Inj\F(\PP)&=\add \sigma_A^{\leq0}\left(\nu_A\PP[-1]\right).
\end{align*}
\item For any indecomposable $P\in\add A$, $\mathfrak{f}_{\PP}(P)$ is either zero or indecomposable, and for indecomposable $P,P'\in\add A$, if $\mathfrak{f}_{\PP}(P)\cong\mathfrak{f}_{\PP}(P')\ne0$ then $P\cong P'$. For any indecomposable $I\in\add \mathbb D A[d-1]$, $\mathfrak{t}_{\PP}(I)$ is either zero or indecomposable, and for indecomposable $I,I'\in\add \mathbb D A[d-1]$, if $\mathfrak{t}_{\PP}(I)\cong\mathfrak{t}_{\PP}(I')\ne0$ then $I\cong I'$.
\item For any indecomposable $P\in\add A$, we have $\mathfrak{f}_{\PP}(P)=0$ if and only if $P\in\add\PP$; $\mathfrak{t}_{\PP}(\nu_A(P)[d-1])=0$ if and only if $P[d]\in\add\PP$. In particular, at least one of $\mathfrak{f}_{\PP}(P)$ and $\mathfrak{t}_{\PP}(\nu_A(P)[d-1])$ is not zero.
\item Let $P_1,\ldots,P_n$ form a complete set of non-isomorphic indecomposable objects in $\add A$. Then
\[|\sigma_A^{\geq-d+1}\PP|= |\mathfrak{t}_{\PP}(\mathbb D A[d-1])|=n-\#\{i\mid P_i[d]\in\add\PP\},\]
\[|\sigma_A^{\leq0}\left(\nu_A(\PP)[-1]\right)|=|\mathfrak{f}_{\PP}(A)|=n-\#\{i\mid P_i\in\add\PP\}.\]
In particular, for each of $\T(\PP)$ and $\F(\PP)$, the numbers of isomorphism classes of indecomposable Ext-projectives and Ext-injectives are equal.
\end{enumerate}
\end{proposition}

\begin{proof}
We prove the assertions concerning $\T(\PP)$; those concerning $\F(\PP)$ follow dually.

For (1), applying \cite[Proposition~2.12(a),(c)]{WZ} to our setting yields
\[
 \Proj\T(\PP)=\add\sigma_A^{\geq-d+1}(\PP),
\]
and shows that $\T(\PP)$ has enough Ext-projectives. The assertions concerning its Ext-injectives follow from Lemma~\ref{lem:torsion-Ext-inj}.

For (2), let $I,I'\in\add(\mathbb D A[d-1])$. Applying $\Hom_{D(A)}(-,I')$ to the canonical triangle $\mathfrak{t}_{\PP} I\to I\to \mathfrak{f}_{\PP} I\to \mathfrak{t}_{\PP} I[1]$ of $I$ with respect to $(\T(\PP),\F(\PP))$ yields an exact sequence
\begin{equation}\label{eq:local1}
    \Hom_{D(A)}(I,I')\to\Hom_{D(A)}(\mathfrak{t}_{\PP}I, I')\to\Hom_{D(A)}(\mathfrak{f}_{\PP}I,I'[1])=0,
\end{equation}
where the last term vanishes since $I'$ is Ext-injective in $\emod{d}{A}$. Applying the functor $\Hom_{D(A)}(\mathfrak{t}_{\PP}I,-)$ to the canonical triangle $\mathfrak{t}_{\PP} I'\to I'\to \mathfrak{f}_{\PP} I'\to \mathfrak{t}_{\PP} I'[1]$ of $I'$ yields an exact sequence
\[\begin{array}{rl}
    & 0=\Hom_{D(A)}(\mathfrak{t}_{\PP}I,\mathfrak{f}_{\PP}I'[-1])\to \Hom_{D(A)}(\mathfrak{t}_{\PP}I,\mathfrak{t}_{\PP}I')\to \Hom_{D(A)}(\mathfrak{t}_{\PP}I, I')\\
    \to & \Hom_{D(A)}(\mathfrak{t}_{\PP}I,\mathfrak{f}_{\PP}I')=0,
\end{array}\]
where the first and last terms vanish since $(\T(\PP),\F(\PP))$ is an $s$-torsion pair. Combining with the surjective map in \eqref{eq:local1}, this implies that $\mathfrak{t}_{\PP}$ induces a surjective map 
\[\Hom_{D(A)}(I,I')\to \Hom_{D(A)}(\mathfrak{t}_{\PP}I,\mathfrak{t}_{\PP}I')\to0.\]
It follows that, if $I$ is indecomposable, then $\End_{D(A)}(\mathfrak{t}_{\PP}I)$ is a quotient of the local algebra $\End_{D(A)}(I)$ and hence is either zero or local. Thus $\mathfrak{t}_{\PP}I$ is zero or indecomposable. Moreover, if $I$ and $I'$ are indecomposable and $\mathfrak{t}_{\PP}(I)\cong\mathfrak{t}_{\PP}(I')\neq 0$, then there exist morphisms $u:I\to I'$ and $v:I'\to I$ such that $\mathfrak{t}_{\PP}(v)\circ\mathfrak{t}_{\PP}(u)=1_{\mathfrak{t}_{\PP}(I)}$. It follows that $v\circ u$ is not nilpotent. Since $\End_{D(A)}(I)$ is local and finite-dimensional, $v\circ u$ is invertible. Thus $u$ is a section, and the indecomposability of $I'$ implies that $I\cong I'$.

For (3), let $P\in\add A$ be indecomposable. By definition, $\mathfrak{t}_{\PP}(\nu_A(P)[d-1])=0$ if and only if $\nu_A(P)[d-1]\in\F(\PP)$, while by \eqref{eq:Serre}, the latter is equivalent to $\Hom_{D(A)}(P[d],\PP[i])=0$ for all $i>0$. Since $\PP\in\add A*\cdots*\add A[d]$, we also have $\Hom_{D(A)}(\PP,P[d+i])=0$ for all $i>0$.
Hence, this is equivalent to $\PP\oplus P[d]$ being silting, and therefore to
$P[d]\in\add\PP$ by \cite[Theorem~2.18]{AI}. 

If both $\mathfrak{f}_{\PP}(P)$ and $\mathfrak{t}_{\PP}(\nu_A(P)[d-1])$ are zero, then both $P$ and $P[d]$ belong to $\add\PP$, contradicting $\Hom_{D(A)}(\PP,\PP[d])=0$.

For (4), we have $|\PP|=|A|=n$ by
\cite[Corollary~2.28]{AI}. Let $r=\#\{i\mid P_i[d]\in\add\PP\}$.  By the equivalence induced by truncation in \cite[Theorem~2.5]{MP}, the indecomposable direct summands of $\PP$ not belonging to $\add A[d]$ correspond bijectively, up to isomorphism, to the indecomposable direct summands of $\sigma_A^{\geq-d+1}(\PP)$. Consequently, $|\sigma_A^{\geq-d+1}(\PP)|=n-r$. By (2) and (3), the non-zero objects $\mathfrak{t}_{\PP}(\nu_A(P_i)[d-1])$ are indecomposable and pairwise non-isomorphic, and exactly $r$ of these $n$ objects vanish. Thus $|\mathfrak{t}_{\PP}(\mathbb D A[d-1])|=n-r$, as required.
\end{proof}

We next recall some facts about almost split triangles. Let $\mathcal{E}\subseteq\emod{d}{A}$ be an extension-closed additive subcategory closed under direct summands. A non-split triangle
\[
 X\xrightarrow{f}M\xrightarrow{g}Z\to X[1]
\]
in $\mathcal{E}$ is called \emph{almost split in $\mathcal{E}$} starting at an indecomposable object $X$ and ending at an indecomposable object $Z$ if, for every $U\in\mathcal{E}$, every morphism $X\to U$ which is not a section factors through $f$, and every morphism $U\to Z$ which is not a retraction factors through $g$. We say that $\mathcal{E}$ \emph{has Auslander-Reiten triangles} if for any indecomposable object $Z$ which is not Ext-projective (resp. Ext-injective) in $\mathcal{E}$, there exists an almost split triangle in $\mathcal{E}$ ending (resp. starting) at it.

By \cite[Theorem~1.2 and Corollary~1.3]{MP}, $\emod{d}{A}$ has Auslander-Reiten triangles. More precisely, for any indecomposable object $Z$ of $\emod{d}{A}$ which is not Ext-projective, there exists an almost split triangle in $\emod{d}{A}$:
\[\tau^A_{[d]}(Z)\to Y\to Z\to \tau^A_{[d]}(Z)[1].\]
For any indecomposable object $X$ of $\emod{d}{A}$ which is not Ext-injective, there exists an almost split triangle in $\emod{d}{A}$:
\[X\to W\to (\tau^A_{[d]})^{-1}X \to X[1].\]
Here, for any object $Z\in\emod{d}{A}$, the object $\tau^A_{[d]}(Z)$ is defined up to isomorphism as follows: by \cite[Theorem~2.5]{MP}, there exists a $(d+1)$-term object $\RR\in\per A$ having no non-zero direct summand in $\add A[d]$ such that $Z\cong\sigma_A^{\geq -d+1}\RR$. Then
\begin{equation}\label{eq:ARd}
    \tau^A_{[d]}(Z):=\sigma_A^{\leq 0}(\nu_A\RR[-1]).
\end{equation}
The object $(\tau^A_{[d]})^{-1}X$ for $X\in\emod{d}{A}$ is defined dually.

The following is the analogue of \cite[Lemma~0.1]{S} in extended module categories.

\begin{lemma}\label{lem:smalo}
Let $(\T,\F)$ be an $s$-torsion pair in
$\emod{d}{A}$.
\begin{enumerate}
\item An indecomposable object $X$ in $\F$ is Ext-injective in $\F$ if and only if either $X\in\Inj\emod{d}{A}$, or $X\notin\Inj\emod{d}{A}$ and $(\tau^A_{[d]})^{-1}X\in\T$.
\item An indecomposable object $Z$ in $\T$ is Ext-projective in $\T$ if and only if either $Z\in\Proj\emod{d}{A}$, or $Z\notin\Proj\emod{d}{A}$ and $\tau^A_{[d]}(Z)\in\F$.
\end{enumerate}
In particular, if $X\to Y\to Z\to X[1]$ is an almost split triangle in $\emod{d}{A}$, with $X\in\F$ and $Z\in\T$, then $X\in\Inj\F$ and $Z\in\Proj\T$.
\end{lemma}

\begin{proof}
We prove (1); assertion (2) is dual. The case $X\in\Inj\emod{d}{A}$ is trivial. Hence, we assume $X\notin\Inj\emod{d}{A}$. Let $Z=(\tau^A_{[d]})^{-1}X$, and let
\[
 X\xrightarrow{f}W\xrightarrow{g} Z \xrightarrow{\delta}X[1]
\]
be the almost split triangle in $\emod{d}{A}$ starting at $X$.

The ``if'' part: suppose $Z\in\T$. For any $\varepsilon\in\Hom_{D(A)}(F,X[1])$ with $F\in\F$, if $\varepsilon\neq 0$, complete it to a triangle $X\xrightarrow{a}E\to F\xrightarrow{\varepsilon}X[1]$ in $\emod{d}{A}$. Since $a$ is not a section, there exists a morphism $b\colon W\to E$ such that $a=b\circ f$. Consequently, there exists a morphism $c\colon Z\to F$ such that the following diagram commutes:
\[\xymatrix{
  X\ar[r]^{f}\ar@{=}[d]&W\ar[r]^{g}\ar[d]^{b}& Z \ar[r]^{\delta}\ar@{-->}[d]^{c}&X[1]\ar@{=}[d]\\
  X\ar[r]_{a}&E\ar[r]&F\ar[r]_{\varepsilon}&X[1]
}\]
Since $Z\in\T$ and $F\in\F$, we have $c=0$, which implies $\delta=\varepsilon\circ c=0$, a contradiction. Thus, $\varepsilon= 0$. Hence,  $X\in\Inj\F$.

The ``only if'' part: suppose that $X\in\Inj\F$ and $Z\notin\T$. Then, in the canonical triangle of $Z$ with respect to $(\T,\F)$
\[
 \mathfrak{t}Z\xrightarrow{\alpha}Z\xrightarrow{\beta}\mathfrak{f}Z
 \to(\mathfrak{t}Z)[1],
\]
the morphism $\alpha$ is not a retraction. It therefore factors through $g$,
which implies that $\delta\circ\alpha=0$. Hence, $\delta$ factors through $\beta$. However, since $\mathfrak{f}Z\in\F$ and $X\in\Inj\F$, we have $\Hom_{D(A)}(\mathfrak{f}Z,X[1])=0$. Therefore, we have $\delta=0$, a contradiction.

The final assertion follows from (1) and (2).
\end{proof}

%==================================================%
\section{Higher Brenner-Butler theorem}\label{sec:higherBB}
%==================================================%

Let $A$ be a proper $d$-truncated non-positive dg algebra. Let $\PP$ be a $(d+1)$-term silting object in $\per A$. In this section, we prove a higher Brenner-Butler theorem for the $s$-torsion pair $(\T(\PP),\F(\PP))$ in $\emod{d}{A}$.

Fix a $K$-projective resolution of $\PP$, which is still denoted by $\PP$. Let $\huaEnd_A(\PP)$ be the dg endomorphism algebra of $\PP$. It is connective since $\PP$ is silting, and is proper since $A$ is proper and $\PP\in\per A$. Let
\[\Gamma=\sigma^{\leq 0}(\huaEnd_A(\PP))\]
be the non-positive truncation of $\huaEnd_A(\PP)$. Thus the inclusion
\[
 j_{\PP}\colon\Gamma\rightarrow\huaEnd_A(\PP)
\]
is a quasi-isomorphism of dg algebras. This inclusion gives $\PP$ a structure of dg left $\Gamma$-module, and the functor
\[
 \RHom_A({}_{\Gamma}\PP,-)\colon D(A)\to D(\Gamma)
\]
is a triangle equivalence, whose quasi-inverse is given by
\[
 -\Lotimes_{\Gamma}\PP\colon D(\Gamma)\to D(A).
\]
See \cite[Lemma~6.1(a)]{Ke1} and \cite[Lemma~3.10]{Ke2}. These restrict to the equivalences
\[\per A\simeq\per \Gamma\text{ and }D_{\mathrm{fd}}(A)\simeq D_{\mathrm{fd}}(\Gamma).\]

Recall that the $d$-extended heart of the $t$-structure induced by $\PP$ is
\[
 \eheart{d}{D_{\mathrm{fd}}(\PP)}
 =\{X\in D_{\mathrm{fd}}(A)\mid
     \Hom_{D(A)}(\PP,X[i])=0\text{ for }i\notin[-d+1,0]\}.
\]
Then the above mutually inverse functors $\RHom_A({}_{\Gamma}\PP,-)$ and $-\Lotimes_{\Gamma}\PP$ restrict to equivalences between $\eheart{d}{D_{\mathrm{fd}}(\PP)}$ and $\emod{d}{\Gamma}$.

\begin{definition}
  The $d$-truncated dg endomorphism algebra of $\PP$ is defined to be
  \[\eend{d}{A}(\PP):=\sigma^{\geq -d+1}(\Gamma).\]
\end{definition}

The $d$-truncated dg endomorphism algebra $\eend{d}{A}(\PP)$ is a proper $d$-truncated non-positive dg algebra. In what follows, we write $B$ for $\eend{d}{A}(\PP)$. Up to quasi-isomorphism, $\Gamma$ and $B$ do not depend on the choice of the $K$-projective resolution of $\PP$.

We write
\[
    p=p_{\PP}\colon\Gamma\twoheadrightarrow B
\]
for the canonical quotient morphism. Let $p_*\colon D(B)\to D(\Gamma)$ be the restriction functor along $p$, and $p^*:=-\Lotimes_{\Gamma}B\colon D(\Gamma)\to D(B)$ be its left adjoint. Composing the functors $p_*$ and $-\Lotimes_{\Gamma}\PP$, we obtain a triangle functor
\begin{equation}\label{eq:real}
  \begin{array}{rcl}
  \real_{\PP}\colon D(B)& \to & D(A)\\
  N &\mapsto &
      p_*(N)\overset{\mathrm L}{\otimes}_\Gamma\PP,
\end{array}
\end{equation}
which restricts to a triangle functor $D_{\mathrm{fd}}(B)\to D_{\mathrm{fd}}(A)$. By the adjunction $p^*\dashv p_*$ and the Hom-tensor adjunction, we have the adjunction
\[
p^*\circ\RHom_A({}_\Gamma\PP,-)\dashv\real_{\PP}.
\]

One can realize the $d$-extended module category $\emod{d}{B}$ via the $d$-extended heart $\eheart{d}{D_{\mathrm{fd}}(\PP)}$ as follows.

\begin{theorem}\label{thm:real-d-heart}
The functor
\[
\begin{array}{rcl}
  \Theta_{\PP}\colon\eheart{d}{D_{\mathrm{fd}}(\PP)}& \to & \emod{d}{B}\\
  M &\mapsto &
      \sigma^{\leq 0}_{B}\RHom_A({}_B B\Lotimes_\Gamma\PP,M)
\end{array}
\]
is an equivalence, whose quasi-inverse is given by $\real_{\PP}$. The equivalence preserves and reflects triangles whose first three terms belong to the respective extended hearts. Moreover, for any $U,V\in\emod{d}{B}$, the functor
$\real_{\PP}$ induces natural isomorphisms
\begin{equation}\label{eq:realization-Hom}
 \Hom_{D(B)}(U,V[i])
 \cong
 \Hom_{D(A)}
   \left(\real_{\PP}(U),\real_{\PP}(V)[i]\right),
 \qquad i\leq1.
\end{equation}
\end{theorem}

\begin{proof}
  This follows directly from Proposition~\ref{prop:trun-equiv} and the equivalence between the categories $\eheart{d}{D_{\mathrm{fd}}(\PP)}$ and $\emod{d}{\Gamma}$, noting that there exists a natural isomorphism 
  \[\RHom_A({}_B B\Lotimes_{\Gamma}\PP,M)\cong \RHom_{\Gamma}\left({}_B B,\RHom_A({}_{\Gamma}\PP,M)\right).\]
\end{proof}

We consider the following full subcategories of $\emod{d}{B}$:
\[\X(\PP):=\Theta_{\PP}(\F(\PP)[d])\quad\text{and}\quad
 \Y(\PP):=\Theta_{\PP}(\T(\PP)).\]

\begin{corollary}\label{cor:higherBB}
The pair $(\X(\PP),\Y(\PP))$ is an $s$-torsion
pair in $\emod{d}{B}$.  There are equivalences $\T(\PP)\simeq\Y(\PP)$ and $\F(\PP)\simeq\X(\PP)$.
\end{corollary}

\begin{proof}
  Since $(\T(\PP),\F(\PP))$ is an $s$-torsion pair in $\emod{d}{A}$,  by \cite[Theorem~1.12]{Z}, the pair $(\F(\PP)[d],\T(\PP))$ is an $s$-torsion pair in the $d$-extended heart $\eheart{d}{D_{\mathrm{fd}}(\PP)}$. Thus, the required assertion follows from Theorem~\ref{thm:real-d-heart}.
\end{proof}

Let $\widetilde{\QQ}=\RHom_A({}_\Gamma\PP,{}_A A[d])\in D(A^{\mathrm{op}}\otimes_{\k}\Gamma)$. We define
\begin{equation}\label{eq:defQ}
  \QQ=p^*\left(\RHom_A({}_\Gamma\PP,A[d])\right)\in D(A^{\mathrm{op}}\otimes_{\k}B).
\end{equation}
Forgetting the left $A$-module structure, we have $\QQ\in D(B)$.

\begin{proposition}\label{prop:Qsilt}
  The object $\QQ$ is a $(d+1)$-term silting object in $\per B$.
\end{proposition}

\begin{proof}
  Since $A[d]$ is a silting object in $\per A$, and by Lemma~\ref{lem:AinP}, $A[d]\in\add\PP*\add\PP[1]*\cdots*\add\PP[d]$, applying the triangle equivalence $\RHom_A({}_\Gamma\PP,-)$, we obtain that $\widetilde{\QQ}$ is a silting object in $\per \Gamma$, and
  \begin{equation}\label{eq:tildeQ}
    \widetilde{\QQ}\in\add\Gamma*\add\Gamma[1]*\cdots*\add\Gamma[d]\subseteq D(\Gamma)^{\leq 0}.
  \end{equation}
  Applying the functor $p^*$ to \eqref{eq:tildeQ} gives $\QQ\in\add B*\add B[1]*\cdots*\add B[d]$. Hence $\QQ$ is a $(d+1)$-term object in $\per B$.

  Applying the functor $\widetilde{\QQ}\Lotimes_{\Gamma}-$ to the canonical exact sequence of dg $\Gamma$-$\Gamma$-bimodules $0\to \ker p\to\Gamma\xrightarrow{p}B\to 0$ gives a triangle in $D(\Gamma)$:
  \begin{equation}\label{eq:triZ}
    Z\to \widetilde{\QQ}\xrightarrow{\eta_{\widetilde{\QQ}}} p_*\QQ\to Z[1],
  \end{equation}
  where $\eta$ is the unit of the adjunction $p^*\dashv p_*$ (cf. \eqref{eq:tensor-unit}), and $Z=\widetilde{\QQ}\Lotimes_{\Gamma}\ker p$. Applying $-\Lotimes_\Gamma\ker p$ to \eqref{eq:tildeQ}, we obtain
\[
 Z\in\add(\ker p)*\add((\ker p)[1])*
 \cdots*\add((\ker p)[d])
 \subseteq D(\Gamma)^{\leq-d}.
\] 
Applying the functor $\Hom_{D(\Gamma)}(\widetilde{\QQ},-)$ to triangle~\eqref{eq:triZ} gives a long exact sequence
  \[
  \cdots\to\Hom_{D(\Gamma)}(\widetilde{\QQ},\widetilde{\QQ}[i])\to\Hom_{D(\Gamma)}(\widetilde{\QQ},p_*\QQ[i])\to\Hom_{D(\Gamma)}(\widetilde{\QQ},Z[i+1])\to\cdots.
  \]
  When $i\geq 1$, the first term is zero because $\widetilde{\QQ}$ is silting, and the last term is zero because $Z\in D(\Gamma)^{\leq-d}$ and $\widetilde{\QQ}\in\add\Gamma*\add\Gamma[1]*\cdots*\add\Gamma[d]$. Hence we have
  \[\Hom_{D(B)}(\QQ,\QQ[i])\cong\Hom_{D(\Gamma)}(\widetilde{\QQ},p_*\QQ[i])=0\quad\text{for all }i>0.\]
  Finally, since $\widetilde{\QQ}$ is silting in $\per \Gamma$, we have $\Gamma\in\thick(\widetilde{\QQ})$. Since $p^*$ is a triangle functor, we have $B=p^*(\Gamma)\in\thick(\QQ)$. Hence $\QQ$ is silting in $\per B$.
\end{proof}

Using the left $A$-module structure of $\QQ$, we obtain an alternative realization of the functor $\real_{\PP}$.

\begin{lemma}\label{lem:real-via-Q}
There is a natural isomorphism of triangle functors
\begin{equation}\label{eq:RHomcongreal}
 \real_{\PP}\cong\RHom_B({}_A\QQ,-)[d]\colon D(B)\to D(A).
\end{equation}
\end{lemma}

\begin{proof}
Let $N\in D(B)$. Since $\RHom_{A}({}_\Gamma\PP,-)$ and $-\Lotimes_{\Gamma}\PP$ are equivalences, the unit \eqref{eq:unit} at $p_*N$ and the counit
\eqref{eq:counit} at ${}_A(A[d])$ give natural isomorphisms
\[p_*N\cong\RHom_A({}_\Gamma\PP,\real_{\PP}(N))\text{ and }\widetilde{\QQ}\Lotimes_\Gamma\PP\cong{}_A(A[d]).\]
Thus we obtain natural isomorphisms in $D(A)$:
\begin{equation}\label{eq:2cong}
\begin{array}{rl}
     & \RHom_B({}_A\QQ,N)\cong\RHom_\Gamma({}_A\widetilde{\QQ},p_*N)\cong\RHom_\Gamma({}_A\widetilde{\QQ},\RHom_A({}_\Gamma\PP,\real_{\PP}(N)))\\
     \cong & \RHom_A(\widetilde{\QQ}\Lotimes_\Gamma\PP,\real_{\PP}(N)) \cong \RHom_A({}_A(A[d]),\real_{\PP}(N)),
\end{array}  
\end{equation}
where the first is given by the adjunction $p^*\dashv p_*$ and the third by the derived tensor-Hom adjunction. Finally, there exists a canonical natural isomorphism from the last term to $\real_{\PP}(N)[-d]$. This proves the assertion.
\end{proof}

Since the object $\QQ$ is a $(d+1)$-term silting object in $\per B$, it induces an $s$-torsion pair $(\T(\QQ),\F(\QQ))$ in $\emod{d}{B}$, where
\[\T(\QQ)=\{X\in\emod{d}{B}\mid \Hom_{D(B)}(\QQ,X[i])=0\text{ for }i>0\},\]
\[\F(\QQ)=\{X\in\emod{d}{B}\mid \Hom_{D(B)}(\QQ,X[i])=0\text{ for }i\leq 0\}.\]

\begin{theorem}\label{thm:higherBB}
  We have $\X(\PP)=\T(\QQ)$ and $\Y(\PP)=\F(\QQ)$. In particular, the $s$-torsion pair $(\X(\PP),\Y(\PP))$ is functorially finite. Moreover, there exist two pairs of mutually quasi-inverse equivalences: 
  \[
  \begin{tikzcd}[column sep=10em, row sep=2em]
    \T(\PP)\arrow[r, shift left=0.6ex, "{\sigma_B^{\leq0}\RHom_A({}_B B\Lotimes_\Gamma\PP,-)}"] & \F(\QQ) \arrow[l, shift left=0.6ex, "{\real_{\PP}\ \cong\ \RHom_B({}_A\QQ,-)[d]}"]\\
    \F(\PP) \arrow[r, shift left=0.6ex, "{\sigma_B^{\leq0}\RHom_A({}_B B\Lotimes_\Gamma\PP,-[d])}"] & \T(\QQ) \arrow[l, shift left=0.6ex, "{\real_{\PP}[-d]\ \cong\ \RHom_B({}_A\QQ,-)}"].
  \end{tikzcd}
\]
\end{theorem}

\begin{proof}
  Let $N\in\emod{d}{B}$. For any $i\in\mathbb Z$, by~\eqref{eq:2cong}, we have natural isomorphisms
  \[
  \Hom_{D(B)}(\QQ,N[i])
 \cong\Hom_{D(A)}(A[d],\real_{\PP}(N)[i]),\quad i\in\mathbb Z.
  \]
  This implies that $N\in\T(\QQ)$ if and only if $\real_{\PP}(N)\in D(A)^{\leq-d}$. By Theorem~\ref{thm:real-d-heart}, the latter is equivalent to $\real_{\PP}(N)\in D(A)^{\leq-d}\cap \eheart{d}{D_{\mathrm{fd}}(\PP)}=\F(\PP)[d]$. Therefore, $\F(\PP)[d] = \real_{\PP}(\T(\QQ))$. Applying the functor $\Theta_{\PP}$ and using Theorem~\ref{thm:real-d-heart}, we obtain $\X(\PP)=\T(\QQ)$. Similarly, we have $\Y(\PP)=\F(\QQ)$.

  By Theorem~\ref{thm:real-d-heart}, the functors displayed from $\T(\PP)$ to $\F(\QQ)$ and from $\F(\PP)$ to $\T(\QQ)$ are equivalences, whose quasi-inverses are the restrictions of $\real_{\PP}$ and $\real_{\PP}(-)[-d]$, respectively. By Lemma~\ref{lem:real-via-Q}, these quasi-inverses are naturally isomorphic to the restrictions of $\RHom_B({}_A\QQ,-)[d]$ and $\RHom_B({}_A\QQ,-)$, respectively.
\end{proof}

%==================================================%
\section{The iterated endomorphism algebra}\label{sec:algC}
%==================================================%

In the previous section, starting from a proper $d$-truncated non-positive dg algebra $A$ and a $(d+1)$-term silting object $\PP$ in $\per A$, we constructed a proper $d$-truncated non-positive dg algebra $B=\eend{d}{A}(\PP)$ and an object $\QQ\in D(A^{\mathrm{op}}\otimes_\k B)$, which, regarded as an object in $D(B)$, is a $(d+1)$-term silting object in $\per B$. Now, we repeat the process with $\QQ$.

Fix a $K$-projective resolution of $\QQ$ in the category of dg $A$-$B$-bimodules, which is still denoted by $\QQ$. Then $\QQ$ is $K$-projective as a right $B$-module (see Section~\ref{subsec:dg-pre}). Let 
\[\Delta=\sigma^{\leq0}(\huaEnd_B(\QQ))\]
be the non-positive truncation of $\huaEnd_B(\QQ)$, and let 
\[C=\eend{d}{B}(\QQ)=\sigma^{\geq -d+1}(\Delta)\]
be the $d$-truncated dg endomorphism algebra of $\QQ$. Write
\[
 j_{\QQ}\colon\Delta\rightarrow\huaEnd_B(\QQ)
 \ \text{and}\ 
 p_{\QQ}\colon\Delta\rightarrow C
\]
for the canonical inclusion and projection, respectively. Then $C$ is a proper $d$-truncated non-positive dg algebra and $j_{\QQ}$ is a quasi-isomorphism.

The left dg $A$-module structure of $\QQ$ gives a dg algebra homomorphism
\[
 \widehat\vartheta\colon A\to\huaEnd_B(\QQ).
\]
Since $A$ is non-positive, there exists a unique morphism $\alpha\colon A\rightarrow\Delta$ such that $j_{\QQ}\circ\alpha=\widehat\vartheta$. We define
\[\vartheta:=p_{\QQ}\circ\alpha\colon A\rightarrow C.\]
See the following commutative diagram.
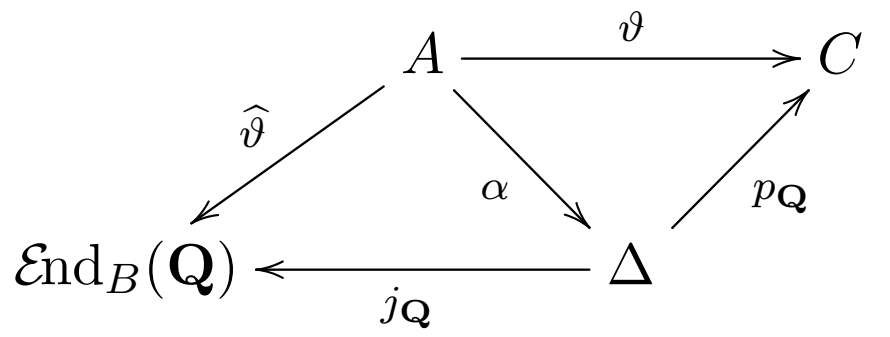
\[\xymatrix{
  &A\ar[rr]^{\vartheta}\ar[dr]_{\alpha}\ar[dl]_{\widehat\vartheta}&&C\\
  \huaEnd_B(\QQ)&&\Delta\ar[ru]_{p_{\QQ}}\ar[ll]^{\quad j_{\QQ}}
}\]

We next describe the restriction functor $\vartheta_*\colon D(C)\to D(A)$. The inclusion $j_{\QQ}$ gives the dg bimodule ${}_\Delta\QQ_B$. Then we have a triangle functor
\[
\begin{array}{rcl}
  \real_{\QQ}\colon D(C)& \to & D(B)\\
  Y &\mapsto &
      (p_{\QQ})_*(Y)\overset{\mathrm L}{\otimes}_\Delta\QQ.
\end{array}
\]

\begin{lemma}\label{lem:res-varthe}
There exists a natural isomorphism of triangle functors
\[
\vartheta_*\cong\real_{\PP}\left(\real_{\QQ}(-)\right)[-d]\colon D(C)\to D(A).
\]
\end{lemma}

\begin{proof}
For any $Y\in D(C)$, since $-\Lotimes_\Delta\QQ\colon D(\Delta)\to D(B)$ is a triangle equivalence, the unit of the derived tensor-Hom adjunction is a natural isomorphism in $D(\Delta)$:
\[(p_{\QQ})_*(Y)\cong\RHom_B
 \left({}_\Delta\QQ,(p_{\QQ})_*(Y)\overset{\mathrm L}{\otimes}_\Delta\QQ\right)=\RHom_B
 \left({}_\Delta\QQ,\real_{\QQ}(Y)\right).\]
Therefore, we obtain the following natural isomorphisms in $D(A)$:
\[\begin{array}{rl}
  & \vartheta_*Y=\alpha_*(p_{\QQ})_*Y\cong\alpha_*\RHom_B\left({}_\Delta\QQ,\real_{\QQ}(Y)\right)\cong\RHom_B\left({}_A\QQ,\real_{\QQ}(Y)\right)\\
  \cong & \real_{\PP}\left(\real_{\QQ}(Y)\right)[-d],
\end{array}
\]
where the second isomorphism follows from $j_{\QQ}\circ\alpha=\widehat\vartheta$, and the last isomorphism follows from Lemma~\ref{lem:real-via-Q}.
\end{proof}

We now characterize when $\vartheta$ is a quasi-isomorphism.

\begin{theorem}\label{thm:equiv}
The following conditions are equivalent:
\begin{enumerate}
\item $\vartheta\colon A\to C$ is a quasi-isomorphism.
\item $p\colon\Gamma\to B$ is a quasi-isomorphism.
\item $\Hom_{\per A}(\PP,\PP[i])=0$ for every $i\leq-d$, i.e., $\Gamma$ is $d$-truncated.
\item $\widehat\vartheta\colon A\to\huaEnd_B(\QQ)$ is a quasi-isomorphism.
\item The restriction functor $\vartheta_*\colon\emod{d}{C}\to\emod{d}{A}$ is an equivalence.
\end{enumerate}
\end{theorem}

\begin{proof}
Applying Proposition~\ref{prop:trun} gives $(1)\Leftrightarrow(5)$. The equivalence $(2)\Leftrightarrow(3)$ follows from
$B=\sigma^{\geq-d+1}(\Gamma)$ and
$H^i(\Gamma)\cong\Hom_{\per A}(\PP,\PP[i])$.
We prove $(2)\Rightarrow(4)\Rightarrow(1)\Rightarrow(2)$.

$(2)\Rightarrow(4)$: By (2),
the functor
$p^*\colon D(\Gamma)\to D(B)$ is an equivalence. Let $\iota\colon\QQ\to p^*\RHom_A({}_{\Gamma}\PP,A[d])$ be the canonical isomorphism in $D(A^{\mathrm{op}}\otimes_\k B)$. For every $n\in\mathbb Z$ and every closed $a\in A^n$, let $\lambda_a\colon A[d]\to A[d+n]$ be the morphism induced by the left action of $a$ on $A[d]$. Then we have the following commutative diagram in $D(B)$:
\[
\xymatrix{
  \QQ\ar[r]^{\iota\qquad\qquad}\ar[d]_{\widehat\vartheta(a)}&p^*\RHom_A({}_{\Gamma}\PP,A[d])\ar[d]^{p^*\RHom_A({}_{\Gamma}\PP,\lambda_a)}\\
  \QQ[n]\ar[r]_{\iota[n]\qquad\qquad}&p^*\RHom_A({}_{\Gamma}\PP,A[d+n])
}
\]
Under the identifications 
\[H^n(A)\cong\Hom_{D(A)}(A[d],A[d+n])\text{ and }H^n\huaEnd_B(\QQ)\cong\Hom_{D(B)}(\QQ,\QQ[n]),\]
the map $H^n(\widehat\vartheta)$ is induced by $p^*\RHom_A({}_\Gamma\PP,-)$. Since this functor is an equivalence, the map $H^n(\widehat\vartheta)$ is an isomorphism for every $n\in\mathbb Z$. Thus $\widehat\vartheta$ is a quasi-isomorphism.

$(4)\Rightarrow(1)$: Since $j_{\QQ}$ is a quasi-isomorphism and
$j_{\QQ}\circ\alpha=\widehat\vartheta$, condition~(4)
implies that $\alpha$ is a quasi-isomorphism.
It follows that $\Delta$ is $d$-truncated. Thus $p_{\QQ}$ is a quasi-isomorphism, and so is
$\vartheta$.

$(1)\Rightarrow(2)$: By (1), $\vartheta_*\colon D(C)\to D(A)$ is an equivalence. Therefore, by Lemma~\ref{lem:res-varthe}, $\real_{\PP}\colon D(B)\to D(A)$ is dense. By \eqref{eq:real}, since $-\overset{\mathrm L}{\otimes}_\Gamma\PP$ is an equivalence, the functor $p_*\colon D(B)\to D(\Gamma)$ is also dense. Hence, there exist $X\in D(B)$ and an isomorphism $u\colon\Gamma\xrightarrow{\cong}p_*X$ in $D(\Gamma)$. By \eqref{eq:adjun-exp}, there exists an isomorphism:
\[
\begin{array}{ccc}
  \Hom_{D(B)}(p^*\Gamma,X) & \xrightarrow{\cong} & \Hom_{D(\Gamma)}(\Gamma,p_*X)\\
  h & \mapsto & p_*(h)\circ \eta_\Gamma
\end{array}
\]
where $\eta$ is the unit of the adjunction $p^*\dashv p_*$. It follows that the morphism 
\[-\circ \eta_\Gamma\colon \Hom_{D(\Gamma)}(p_*p^*\Gamma,p_*X)\to\Hom_{D(\Gamma)}(\Gamma,p_*X)\]
is surjective. Thus, the isomorphism $u$ factors through $p$, which implies that $\eta_\Gamma\colon\Gamma\to p_*p^*\Gamma\cong p_*B$ is a section. Thus $\Gamma$ is a direct summand of $p_*B$.
Therefore, due to $p_*B\in D(\Gamma)^{\geq-d+1}$,
we have $\Gamma\in D(\Gamma)^{\geq-d+1}$. Hence $p\colon\Gamma\to B$ is a quasi-isomorphism.
\end{proof}

Applying Corollary~\ref{cor:higherBB} to $\QQ$, we obtain an $s$-torsion pair $(\X(\QQ),\Y(\QQ))$ in $\emod{d}{C}$, where
\[\X(\QQ)=\Theta_{\QQ}(\F(\QQ)[d])\ \text{and}\ \Y(\QQ)=\Theta_{\QQ}(\T(\QQ)).\]

\begin{theorem}\label{thm:torsion-res}
The functor $\vartheta_*\colon\emod{d}{C}\to\emod{d}{A}$ induces equivalences
\[
 \X(\QQ)
 \xrightarrow{\simeq}
 \T(\PP)
 \quad\text{and}\quad
 \Y(\QQ)
 \xrightarrow{\simeq}
 \F(\PP).
\]
\end{theorem}

\begin{proof}
Applying Theorem~\ref{thm:higherBB} to $\PP$ and $\QQ$, we obtain equivalences
\[
\begin{gathered}
 \X(\QQ)\xrightarrow{\real_{\QQ}(-)[-d]}\F(\QQ)=\Y(\PP)\xrightarrow{\real_{\PP}}\T(\PP),\\
 \Y(\QQ)\xrightarrow{\real_{\QQ}}\T(\QQ)=\X(\PP)\xrightarrow{\real_{\PP}(-)[-d]}\F(\PP).
\end{gathered}
\]
Both compositions are naturally isomorphic to
$\real_{\PP}\left(\real_{\QQ}(-)\right)[-d]$ and hence to $\vartheta_*$ by Lemma~\ref{lem:res-varthe}.
\end{proof}

This result recovers \cite[Theorem 4.4]{BZ} for the case $d=1$. Note that, in that case, $\vartheta_*\colon\emod{d}{C}\to \emod{d}{A}$ is always fully faithful. However, this is not true for general $d$. See the following example.

\begin{example}\label{ex:A3-worked}
Let $d=2$, and $A=\k Q$, where $Q$ is the quiver 
\[
1\xrightarrow{\alpha}2\xrightarrow{\beta}3
\]
with $|\alpha|=|\beta|=0$. Let $P_i^A$ (resp. $I_i^A$ and $S_i^A$) denote the indecomposable projective (resp. indecomposable injective and simple) right $A$-module corresponding to vertex $i$. We use analogous notation for algebras $B$ and $C$ below. Let
\[
\PP=P_3^A\oplus P_2^A[1]\oplus P_1^A[2].
\]
It is straightforward to check that $\PP$ is a 3-term silting object in $\per A$, and the induced $s$-torsion pair $(\T(\PP),\F(\PP))$ in $\emod{2}{A}$ is given by
\[
\begin{aligned}
 \T(\PP)
   &=\add\left(P_3^A\oplus P_3^A[1]\oplus P_2^A[1]\oplus S_2^A[1]\right),\\
 \F(\PP)
   &=\add\left(S_2^A\oplus I_2^A\oplus I_1^A\oplus I_1^A[1]\right).
\end{aligned}
\]
See Figure~\ref{fig:1}, where red and blue shaded regions indicate the torsion and torsion-free classes of the displayed $s$-torsion pair, respectively.
\begin{figure}[htbp]
\begin{center}
\begin{tikzpicture}
   \begin{scope}[on grid,node distance=2cm]
    \node (A) {$P_2^A$};
    \node (B)[right = 2cm of A] {$I_2^A$};
    \node (C)[above right = 0.9cm and 1cm of A] {$P_1^A$};
    \node (D)[right = 2cm of C] {$P_3^A[1]$};
    \node (E)[right = 2cm of B] {$P_2^A[1]$};
    \node (F)[right = 2cm of D] {$S_2^A[1]$};
    \node (G)[below left = 0.9cm and 1cm of A] {$P_3^A$};
    \node (H)[right = 2cm of G] {$S_2^A$};
    \node (I)[right = 2cm of H] {$I_1^A$};
    \node (J)[right = 2cm of I] {$P_1^A[1]$};
    \node (K)[right = 2cm of E] {$I_2^A[1]$};
    \node (L)[right = 2cm of F] {$I_1^A[1]$};
    \path   (A) edge[->]    (C)
        (C) edge[->]   (B)
        (B) edge[->]  (D);
      \path  (D) edge[->]  (E)
        (E) edge[->]  (F)
        (G) edge[->] (A)
        (A) edge[->] (H)
        (H) edge[->] (B)
        (B) edge[->] (I)
        (I) edge[->] (E)
        (E) edge[->] (J)
        (J) edge[->] (K)
        (F) edge[->] (K)
        (K) edge[->] (L)
        (A) edge[-, dashed] (B)
         (B) edge[-, dashed] (E)
         (C) edge[-, dashed] (D)
         (D) edge[-, dashed] (F)
         (F) edge[-, dashed] (L)
         (E) edge[-, dashed] (K)
         (G) edge[-, dashed] (H)
         (H) edge[-, dashed] (I)
         (I) edge[-, dashed] (J);
    \end{scope}
     \begin{scope}[on background layer]
    \draw[red, thick, fill=red!10, draw=red]
        (G) circle (14pt);
\end{scope}
\begin{scope}[on background layer]
\draw[blue, thick, fill=blue!10, draw=blue]
        (L) circle (14pt);
\end{scope}
\begin{scope}[on background layer]
    \draw[
        blue,
        thick,
        rounded corners=12pt,
        fill=blue!10,
         fill opacity=0.8
    ]
        ([shift={(-12pt,-2pt)}]H.south west) --
        ([shift={(0pt,16pt)}]B.north) --
        ([shift={(12pt,-2pt)}]I.south east) --
        cycle;
\end{scope}
\begin{scope}[on background layer]
  \draw[
    red,
    thick,
    rounded corners=12pt,
    fill=red!10,
    fill opacity=0.8
]
    ([shift={(-12pt, 2pt)}]D.north west) --
    ([shift={(0pt, -16pt)}]E.south) --
    ([shift={(12pt, 2pt)}]F.north east) --
    cycle;
\end{scope}
   \end{tikzpicture}
\end{center}
\caption{$(\T(\PP),\F(\PP))$ in $\emod{2}{A}$.}\label{fig:1}
\end{figure}

Label the primitive idempotents of $\Gamma=\huaEnd_A(\PP)$ so
that the vertices $1,2,3$ correspond to
$P_1^A[2],P_2^A[1],P_3^A$, respectively. Then $\Gamma$ has zero
differential and
\[
 \Gamma\cong
 \k\left(
  1\xrightarrow{a}2
   \xrightarrow{b}3
 \right)
\]
as a graded dg algebra, where $|a|=|b|=-1$. The path $ab$ has degree $-2$, so
\[
 B=\sigma^{\geq -1}(\Gamma)
 \cong
 \k\left(
  1\xrightarrow{a}2
   \xrightarrow{b}3
 \right)/(ab),
\]
again with zero differential. The induced 3-term silting object in $\per B$ is
\[\QQ\cong P_1^B\oplus P_2^B[1]\oplus P_3^B[2],\]
and the induced $s$-torsion pair $(\T(\QQ),\F(\QQ))=(\X(\PP),\Y(\PP))$ is given by
\[
\begin{aligned}
  \T(\QQ)
   &=\add\left(S_2^B[1]\oplus S_1^B\oplus S_1^B[1]\oplus P_1^B\right),\\
 \F(\QQ)
   &=\add\left(P_3^B\oplus P_3^B[1]\oplus S_2^B\oplus P_2^B\right).
\end{aligned}
\]
See Figure~\ref{fig:2}.
\begin{figure}[htbp]
\begin{centering}
\begin{tikzpicture}
   \begin{scope}[on grid,node distance=2cm]
    \node (A) {$P_3^B$};
    \node (B)[right = 2cm of A] {$P_3^B[1]$};
    \node (C)[right = 2cm of B] {$S_2^B$};
    \node (D)[right = 2cm of C] {$S_2^B[1]$};
    \node (E)[right = 2cm of D] {$S_1^B$};
    \node (F)[right = 2cm of E] {$S_1^B[1]$};
    \node (G)[above right = 0.9cm and 1cm of B] {$P_2^B$};
    \node (H)[above right = 0.9cm and 1cm of D] {$P_1^B$};

    \path
        (B) edge[->] (G)
        (G) edge[->] (C)
        (D) edge[->] (H)
        (H) edge[->] (E)
        (A) edge[-, dashed] (B)
        (B) edge[-, dashed] (C)
        (C) edge[-, dashed] (D)
        (D) edge[-, dashed] (E)
        (E) edge[-, dashed] (F);
\end{scope}
   \begin{scope}[on background layer]
    \draw[
        red,
        thick,
        rounded corners=15pt,
        fill=red!10,
        fill opacity=1
    ]
        ([shift={(-10pt,10pt)}]D.north west) --
        ([shift={(-10pt,-10pt)}]D.south west) --
        ([shift={(10pt,-10pt)}]F.south east) --
        ([shift={(10pt,10pt)}]F.north east) --
        ([shift={(0pt,20pt)}]E.north) --
        ([shift={(0pt,30pt)}]H.north) --
        cycle;
   \end{scope}
   \begin{scope}[on background layer]
    \draw[
        blue,
        thick,
        rounded corners=15pt,
        fill=blue!10,
        fill opacity=0.8
    ]
        ([shift={(-10pt,10pt)}]A.north west) --
        ([shift={(-10pt,-10pt)}]A.south west) --
        ([shift={(10pt,-10pt)}]C.south east) --
        ([shift={(10pt,10pt)}]C.north east) --
        ([shift={(0pt,30pt)}]G.north) --
        ([shift={(0pt,20pt)}]B.north) --
        cycle;
   \end{scope}
\end{tikzpicture}
\end{centering}
\caption{$(\T(\QQ),\F(\QQ))$ in $\emod{2}{B}$.}\label{fig:2}
\end{figure}

The $d$-truncated dg endomorphism algebra $C$ of $\QQ$ is quasi-isomorphic to
\[\k(1\xrightarrow{\alpha}2\xrightarrow{\beta}3)/(\alpha\beta),
\]
where $|\alpha|=|\beta|=0$. Under this quasi-isomorphism, the morphism $\vartheta\colon A\to C$ is the canonical quotient. Moreover, the induced $s$-torsion pair $(\X(\QQ),\Y(\QQ))$ is given by
\[
\begin{aligned}
 \X(\QQ)
   &=\add\left(S_3^C\oplus S_3^C[1]\oplus P_2^C[1]
                  \oplus S_2^C[1]\right),\\
 \Y(\QQ)
   &=\add\left(S_1^C\oplus S_1^C[1]\oplus S_2^C\oplus P_1^C\right).
\end{aligned}
\]
See Figure~\ref{fig:3}, where
$X=(P_2^C\xrightarrow{\alpha}P_1^C)$ is concentrated in degrees
$-1,0$.
\begin{figure}[htbp]
\centering
\begin{tikzpicture}
    \begin{scope}[on grid,node distance=2cm]
        \node (A) {$P_3^C$};
        \node (B)[right = 2cm of A] {$S_2^C$};
        \node (C)[right = 2cm of B] {$X$};
        \node (D)[right = 2cm of C] {$S_2^C[1]$};
        \node (E)[right = 2cm of D] {$S_1^C[1]$};

        \node (F)[above right = 0.9cm and 1cm of B] {$P_1^C$};
        \node (G)[right = 2cm of F] {$P_2^C[1]$};

        \node (H)[below right = 0.9cm and 1cm of A] {$P_2^C$};
        \node (I)[right = 2cm of H] {$P_3^C[1]$};
        \node (J)[right = 2cm of I] {$S_1^C$};
        \node (K)[right = 2cm of J] {$P_1^C[1]$};

        \path
            (A) edge[->] (H)
            (H) edge[->] (B)
            (B) edge[->] (I)
            (B) edge[->] (F)
            (F) edge[->] (C)
            (I) edge[->] (C)
            (C) edge[->] (J)
            (C) edge[->] (G)
            (J) edge[->] (D)
            (G) edge[->] (D)
            (D) edge[->] (K)
            (K) edge[->] (E)
            (A) edge[-, dashed] (B)
            (B) edge[-, dashed] (C)
            (C) edge[-, dashed] (D)
            (D) edge[-, dashed] (E)
            (F) edge[-, dashed] (G)
            (H) edge[-, dashed] (I)
            (I) edge[-, dashed] (J)
            (J) edge[-, dashed] (K);
    \end{scope}

    \begin{scope}[on background layer]
        \draw[draw=red, thick, fill=red!10]
            (A) circle (16pt);
        \draw[draw=red, thick, fill=red!10]
            (I) circle (16pt);
        \draw[draw=blue, thick, fill=blue!10]
            (J) circle (16pt);
        \draw[draw=blue, thick, fill=blue!10]
            (E) circle (16pt);
        \draw[
            draw=blue,
            thick,
            fill=blue!10,
            fill opacity=1,
            rotate around={38:($(B)!0.5!(F)+(0pt,2pt)$)}
        ]
        ($(B)!0.5!(F)+(0pt,2pt)$) ellipse [x radius=1.15cm, y radius=0.72cm];
        \draw[
            draw=red,
            thick,
            fill=red!10,
            fill opacity=1,
            rotate around={142:($(G)!0.5!(D)+(0pt,2pt)$)}
        ]
        ($(G)!0.5!(D)+(0pt,2pt)$) ellipse [x radius=1.15cm, y radius=0.72cm];
    \end{scope}
\end{tikzpicture}

\caption{$(\X(\QQ),\Y(\QQ))$ in $\emod{2}{C}$.}
\label{fig:3}
\end{figure}

Although the functor $\vartheta_*$ restricts to equivalences between $\X(\QQ)$ and $\T(\PP)$ and between $\Y(\QQ)$ and $\F(\PP)$ by Theorem~\ref{thm:torsion-res}, it does not restrict to a fully faithful functor from $\emod{2}{C}$ to $\emod{2}{A}$ in this example. This is because the image of $X$ under $\vartheta_*$ is decomposable: $\vartheta_*X\cong P_3^A[1]\oplus I_1^A$. To see this, consider the canonical triangle of $X$ with respect to the $s$-torsion pair $(\X(\QQ),\Y(\QQ))$:
\[P_3^C[1]\to X\to S_1^C.\]
Its image under $\vartheta_*$ is $P_3^A[1]\to \vartheta_*X\to I_1^A=S_1^A$, which is split.
\end{example}

%==================================================%
\section{Auslander-Reiten theory}\label{sec:AR-theory}
%==================================================%

In this section, we apply the higher Brenner-Butler theorem to study the AR theory of $\emod{d}{B}$, where $B=\eend{d}{A}(\PP)$.

For an indecomposable $P\in\add A$, we write
$I_P=\nu_A(P)[d-1]\in\Inj\emod{d}{A}$. Recall that $\mathfrak{t}_{\PP}$ and $\mathfrak{f}_{\PP}$ denote torsion and torsion-free functors respectively, with respect to the $s$-torsion pair $(\T(\PP),\F(\PP))$. 

\begin{lemma}\label{lem:equivalences}
    For any indecomposable $P\in\add A$, at least one of $\Theta_\PP((\mathfrak{f}_{\PP}P)[d])$ and $\Theta_\PP(\mathfrak{t}_{\PP}I_P)$ is not zero. Moreover, any non-zero object among them is indecomposable.
\end{lemma}

\begin{proof}
    By Theorem~\ref{thm:real-d-heart}, $\Theta_\PP$ is an equivalence. Therefore, this lemma follows from Proposition~\ref{prop:Ext-proj-inj}~(2) and (3).
\end{proof}

For an indecomposable $P\in\add A$, we define $\widetilde \QQ_P=\RHom_A({}_\Gamma\PP,P[d])$ and $\QQ_P=p^*(\widetilde \QQ_P)$. They are direct summands of $\widetilde{\QQ}$ and $\QQ$, respectively.

\begin{lemma}\label{lem:QP}
  For any indecomposable $P\in\add A$, $\QQ_P$ is indecomposable, and
  \begin{equation}\label{eq:con-end-cong}
 \sigma_B^{\geq-d+1}(\QQ_P)
 \cong\Theta_{\PP}((\mathfrak{f}_{\PP}P)[d])\text{ and }
 \sigma_B^{\leq0}\left(\nu_B(\QQ_P)[-1]\right)
 \cong\Theta_{\PP}(\mathfrak{t}_{\PP}I_P).
\end{equation}
\end{lemma}

\begin{proof}
  For $N\in\emod{d}{B}$, taking cohomology in \eqref{eq:2cong} and then the direct summands corresponding
to $P$ gives natural isomorphisms
\begin{equation}\label{eq:QPN}
  \Hom_{D(B)}(\QQ_P,N[i])\cong\Hom_{D(A)}(P[d],\real_{\PP}(N)[i]),\quad i\in\mathbb Z.
\end{equation}
Applying $\Hom_{D(A)}(-,\real_{\PP}(N))$ to the canonical triangle $\mathfrak{t}_{\PP}P\to P\to\mathfrak{f}_{\PP}P\to\mathfrak{t}_{\PP}P[1]$ of $P$ with respect to $(\T(\PP),\F(\PP))$ yields an exact sequence
\[
\begin{array}{rl}
  &0=\Hom_{D(A)}\left((\mathfrak{t}_{\PP}P)[d+1],\real_{\PP}(N)\right)\to \Hom_{D(A)}\left((\mathfrak{f}_{\PP}P)[d],\real_{\PP}(N)\right)\\
  \to & \Hom_{D(A)}\left(P[d],\real_{\PP}(N)\right)\to \Hom_{D(A)}\left((\mathfrak{t}_{\PP}P)[d],\real_{\PP}(N)\right)=0,
\end{array}
\]
where the first and last terms vanish since $\mathfrak{t}_{\PP}P\in\T(\PP)\subseteq D_{\mathrm{fd}}(\PP)^{\leq 0}$ and $\real_{\PP}(N)\in \eheart{d}{D_{\mathrm{fd}}(\PP)}\subseteq D_{\mathrm{fd}}(\PP)^{\geq-d+1}$. Consequently, we obtain natural isomorphisms
\[\begin{array}{l}
  \Hom_{D(B)}(\sigma_B^{\geq-d+1}(\QQ_P),N)\cong\Hom_{D(B)}(\QQ_P,N)\overset{\eqref{eq:QPN}}{\cong}
   \Hom_{D(A)}\left(P[d],\real_{\PP}(N)\right)\\
  \cong \Hom_{D(A)}
   \left((\mathfrak{f}_{\PP}P)[d],\real_{\PP}(N)\right)\overset{\eqref{eq:realization-Hom}}{\cong}
   \Hom_{D(B)}
   \left(\Theta_{\PP}((\mathfrak{f}_{\PP}P)[d]),N\right),
\end{array}\]
where the first isomorphism is given by canonical truncation of $\QQ_P\in D(B)^{\leq0}$. Then the Yoneda Lemma shows the first isomorphism in \eqref{eq:con-end-cong}.

Next, applying $\Hom_{D(A)}(\real_{\PP}(N),-)$ to the canonical triangle $\mathfrak{t}_{\PP}I_P\to I_P\to\mathfrak{f}_{\PP}I_P\to\mathfrak{t}_{\PP}I_P[1]$ of $I_P$ with respect to $(\T(\PP),\F(\PP))$ yields an exact sequence
\[\begin{array}{rl}
  &0=\Hom_{D(A)}
   \left(\real_{\PP}(N),\mathfrak{f}_{\PP}I_P[-1]\right)\to \Hom_{D(A)}
   \left(\real_{\PP}(N),\mathfrak{t}_{\PP}I_P\right) \\
   \to & \Hom_{D(A)}\left(\real_{\PP}(N),I_P\right)\to \Hom_{D(A)}
   \left(\real_{\PP}(N),\mathfrak{f}_{\PP}I_P\right)=0,
\end{array}\]
where the first and last terms vanish since $\mathfrak{f}_{\PP}I_P\in\F(\PP)\subseteq D_{\mathrm{fd}}(\PP)^{\geq1}$ and $\real_{\PP}(N)\in\eheart{d}{D_{\mathrm{fd}}(\PP)}\subseteq D_{\mathrm{fd}}(\PP)^{\leq0}$. Consequently, we obtain natural isomorphisms
\[\begin{array}{l}
  \Hom_{D(B)}(N,\sigma_B^{\leq0}\left(\nu_B(\QQ_P)[-1]\right))\cong\Hom_{D(B)}(N,\nu_B(\QQ_P)[-1])\overset{\eqref{eq:Serre}}{\cong}
   \mathbb D\Hom_{D(B)}(\QQ_P,N[1])\\
   \overset{\eqref{eq:QPN}}{\cong}
   \mathbb D\Hom_{D(A)}
   \left(P[d],\real_{\PP}(N)[1]\right)
   \overset{\eqref{eq:Serre}}{\cong}\Hom_{D(A)}
   \left(\real_{\PP}(N)[1],\nu_A(P[d])\right)\\
   \cong
   \Hom_{D(A)}
   \left(\real_{\PP}(N),\mathfrak{t}_{\PP}I_P\right)\overset{\eqref{eq:realization-Hom}}{\cong}
   \Hom_{D(B)}
   \left(N,\Theta_{\PP}(\mathfrak{t}_{\PP}I_P)\right),
\end{array}\]
where the first isomorphism is given by canonical truncation of $\nu_B(\QQ_P)[-1]\in D(B)^{\geq-d+1}$. Then the Yoneda Lemma shows the second isomorphism in \eqref{eq:con-end-cong}.

By Lemma~\ref{lem:equivalences}, the isomorphisms in \eqref{eq:con-end-cong} imply that at least one of $\sigma_B^{\geq-d+1}(\QQ_P)$ and $\sigma_B^{\leq0}\left(\nu_B(\QQ_P)[-1]\right)$ is not zero. It follows that $\QQ_P\neq 0$. Taking the direct summand corresponding to $P$ in
\eqref{eq:triZ}, we obtain a triangle
\[
 Z_P\to\widetilde \QQ_P\to p_*\QQ_P\to Z_P[1],
\]
with $Z_P\in D(\Gamma)^{\leq-d}$. By $\widetilde \QQ_P\in\add\Gamma*\cdots*\add\Gamma[d]$, we have $\Hom_{D(\Gamma)}(\widetilde \QQ_P,Z_P[1])=0$. Applying $\Hom_{D(\Gamma)}(\widetilde \QQ_P,-)$ to the above triangle and using the adjunction $p^*\dashv p_*$ gives a surjective algebra homomorphism
\[
 \End_{D(A)}(P)\cong\End_{D(\Gamma)}(\widetilde \QQ_P)\twoheadrightarrow\End_{D(B)}(\QQ_P).
\]
Since $\End_{D(A)}(P)$ is local and $\QQ_P\ne0$, the object $\QQ_P$ is indecomposable.
\end{proof}

An almost split triangle in $\emod{d}{B}$ is called \emph{connecting} if its left term belongs to $\Y(\PP)$ and its right term belongs to $\X(\PP)$.

\begin{lemma}\label{lem:connecting1}
  If $Y\to E\to X\to Y[1]$ is a connecting triangle, then there exists an indecomposable $P\in\add A$ such that $Y\cong\Theta_{\PP}(\mathfrak{t}_{\PP}I_P)$.
\end{lemma}

\begin{proof}
By Lemma~\ref{lem:smalo}, $Y\in\Inj\Y(\PP)$. By Theorem~\ref{thm:real-d-heart} and Corollary~\ref{cor:higherBB}, there exists an indecomposable $Z\in\Inj\T(\PP)$ such that $Y\cong\Theta_{\PP}(Z)$. By Proposition~\ref{prop:Ext-proj-inj}, $\Inj\T(\PP)=\add\mathfrak{t}_{\PP}(\mathbb D A[d-1])$. Thus, the assertion follows.
\end{proof}

We now compute the Auslander-Reiten translations of the end terms
of connecting triangles.

\begin{lemma}\label{lem:connecting2}
If $\Theta_{\PP}((\mathfrak{f}_{\PP}P)[d])\neq 0$, then
\[
 \tau^B_{[d]}\Theta_{\PP}((\mathfrak{f}_{\PP}P)[d])\cong \Theta_{\PP}(\mathfrak{t}_{\PP}I_P).
\]
If $\Theta_{\PP}(\mathfrak{t}_{\PP}I_P)\neq 0$, then
\[
 (\tau^B_{[d]})^{-1}\Theta_{\PP}(\mathfrak{t}_{\PP}I_P)\cong\Theta_{\PP}((\mathfrak{f}_{\PP}P)[d]).
\]
In particular, if $\Theta_{\PP}((\mathfrak{f}_{\PP}P)[d])=0$, then the object $\Theta_{\PP}(\mathfrak{t}_{\PP}I_P)$ is Ext-injective in $\emod{d}{B}$;
if $\Theta_{\PP}(\mathfrak{t}_{\PP}I_P)=0$, then the object $\Theta_{\PP}((\mathfrak{f}_{\PP}P)[d])$ is Ext-projective in $\emod{d}{B}$.
\end{lemma}

\begin{proof}
  Assume that $\Theta_{\PP}((\mathfrak{f}_{\PP}P)[d])\neq0$. By \eqref{eq:con-end-cong}, $\sigma_B^{\geq-d+1}(\QQ_P)\neq0$. Since $\QQ_P$ is indecomposable by Lemma~\ref{lem:QP}, it has no non-zero direct summand in $\add B[d]$. Thus we have
  \[\tau^B_{[d]}\left(\Theta_{\PP}((\mathfrak{f}_{\PP}P)[d])\right)\overset{\eqref{eq:con-end-cong}}{\cong}\tau^B_{[d]}\left(\sigma_B^{\geq-d+1}(\QQ_P)\right)\overset{\eqref{eq:ARd}}{\cong}\sigma_B^{\leq0}\left(\nu_B(\QQ_P)[-1]\right)\overset{\eqref{eq:con-end-cong}}{\cong}\Theta_{\PP}(\mathfrak{t}_{\PP}I_P).\]
  The other isomorphism can be shown dually.

  If $\Theta_{\PP}(\mathfrak{t}_{\PP}I_P)=0$, then $\Theta_{\PP}((\mathfrak{f}_{\PP}P)[d])\neq 0$. Hence, by the first isomorphism in this lemma, we have $\tau^B_{[d]}\Theta_{\PP}((\mathfrak{f}_{\PP}P)[d])=0$, which implies that  $\Theta_{\PP}((\mathfrak{f}_{\PP}P)[d])$ is Ext-projective in $\emod{d}{B}$. The assertion for the case $\Theta_{\PP}((\mathfrak{f}_{\PP}P)[d])=0$ can be shown dually.
\end{proof}

By Lemmas~\ref{lem:connecting1} and~\ref{lem:connecting2}, together with Proposition~\ref{prop:Ext-proj-inj}~(3), the connecting almost split triangles in $\emod{d}{B}$ are precisely those of the form
\begin{equation}\label{eq:connecting-triangle}
 \Theta_{\PP}(\mathfrak{t}_{\PP}I_P)\xrightarrow{l_P}E_P\xrightarrow{r_P}\Theta_{\PP}((\mathfrak{f}_{\PP}P)[d])
 \xrightarrow{\delta_P}\Theta_{\PP}(\mathfrak{t}_{\PP}I_P)[1],
 \quad P,P[d]\notin\add\PP.
\end{equation}

\begin{proposition}\label{prop:connecting}
Let $P\in\add A$ be indecomposable with
$P,P[d]\notin\add\PP$.
The canonical triangle of the middle term $E_P$ in
\eqref{eq:connecting-triangle} with respect to $(\X(\PP),\Y(\PP))$ is
\[
 \Theta_{\PP}((\operatorname{rad}\mathfrak{f}_{\PP}P)[d])\to E_P \to \Theta_{\PP}(\mathfrak{t}_{\PP}I_P/\operatorname{soc}\mathfrak{t}_{\PP}I_P) \to\Theta_{\PP}((\operatorname{rad}\mathfrak{f}_{\PP}P)[d])[1].
\]
\end{proposition}

\begin{proof}
We adapt the proof of \cite[Corollary~5.3]{BZ} with suitable modifications.

Let $S=\operatorname{top}_{H^0(A)}H^0(P)$, regarded as a dg $A$-module via $A\to H^0(A)$. Since $H^0(P)$ is an indecomposable projective $H^0(A)$-module by \eqref{eq:proj-equiv}, $S$ is a simple $H^0(A)$-module. Since $P,P[d]\notin\add\PP$, by Proposition~\ref{prop:Ext-proj-inj}~(3), both $\mathfrak{f}_{\PP}P$ and $\mathfrak{t}_{\PP}I_P$ are non-zero. We claim that
\[
 \operatorname{top}\mathfrak{f}_{\PP}P\cong S\text{ and } \operatorname{soc}\mathfrak{t}_{\PP}I_P\cong S[d-1].
\]
Consider the canonical triangle of $P$ with respect to $(\T(\PP),\F(\PP))$, after a rotation:
\[\mathfrak{f}_{\PP}P[-1]\to \mathfrak{t}_{\PP}P\to P\to \mathfrak{f}_{\PP}P.\]
Taking $H^0$ gives a surjection $H^0(P)\to H^0(\mathfrak{f}_{\PP}P)$. If $H^0(\mathfrak{f}_{\PP}P)=0$, then $\mathfrak{f}_{\PP}P[-1]\in\emod{d}{A}$. Since $P\in\Proj\emod{d}{A}$, it follows that the above triangle splits. Thus,  $P\in\T(\PP)$ and hence $\mathfrak{f}_{\PP}P=0$, a contradiction. Thus $H^0(\mathfrak{f}_{\PP}P)$ is a non-zero quotient of the indecomposable projective $H^0(P)$, so its top is $S$. Dually, $H^{-d+1}(\mathfrak{t}_{\PP}I_P)$ is a non-zero
submodule of $H^{-d+1}(I_P)$, whose socle is $S$ by
\eqref{eq:Serre}. This proves the claim.

Consider the following two triangles in $\emod{d}{A}$ (see Definition~\ref{def:rad-soc}):
\[
  \operatorname{rad}\mathfrak{f}_{\PP}P\xrightarrow{\alpha}\mathfrak{f}_{\PP}P\xrightarrow{\beta}S\xrightarrow{\gamma}(\operatorname{rad}\mathfrak{f}_{\PP}P)[1],
\]
\[
  S[d-1]\xrightarrow{f}\mathfrak{t}_{\PP}I_P\xrightarrow{g}\mathfrak{t}_{\PP}I_P/\operatorname{soc}\mathfrak{t}_{\PP}I_P\xrightarrow{h} S[d].
\]
Since $S[-1]\in D_{\mathrm{fd}}(A)^{\geq 1}\subseteq D_{\mathrm{fd}}(\PP)^{\geq1}$ and $S[d]\in D_{\mathrm{fd}}(A)^{\leq -d}\subseteq D_{\mathrm{fd}}(\PP)^{\leq0}$ by \eqref{eq:silt-inc}, applying $\Hom_{D(A)}(\PP,-)$ to these two triangles shows that $\operatorname{rad}\mathfrak{f}_{\PP}P\in\F(\PP)$ and $\mathfrak{t}_{\PP}I_P/\operatorname{soc}\mathfrak{t}_{\PP}I_P\in\T(\PP)$. 

Since $\beta\neq 0$, the morphism $\alpha$ is not a retraction. Hence $\Theta_{\PP}(\alpha[d])$ factors through the morphism $r_P$ in \eqref{eq:connecting-triangle}. Consequently, $\delta_P\circ\Theta_{\PP}(\alpha[d])=0$. Applying $\real_{\PP}$, by Theorem~\ref{thm:real-d-heart}, we have $\real_{\PP}(\delta_P)\circ\alpha[d]=0$. Hence, there exists a morphism $\theta\colon S[d-1]\to \mathfrak{t}_{\PP}I_P$ such that $\real_{\PP}(\delta_P)=\theta[1]\circ\beta[d]$; see the following commutative diagram.
\[\xymatrix{
(\operatorname{rad}\mathfrak{f}_{\PP}P)[d]\ar[r]^{\quad\alpha[d]} & (\mathfrak{f}_{\PP}P)[d]\ar[r]^{\quad\beta[d]}\ar[d]_{\real_{\PP}(\delta_P)} & S[d]\ar[r]^{\gamma[d]\qquad\quad}\ar@{-->}[ld]^{\theta[1]} & (\operatorname{rad}\mathfrak{f}_{\PP}P)[d+1]\\
&\mathfrak{t}_{\PP}I_P[1]
}\]
Since $\real_{\PP}(\delta_P)\neq 0$ by \eqref{eq:realization-Hom}, we have $\theta\neq 0$. Then by \eqref{eq:simple-socle-Hom} and $\operatorname{soc}\mathfrak{t}_{\PP}I_P\cong S[d-1]$, there exists an automorphism $u$ of $S$ such that $\theta=f\circ (u[d-1])$. Consequently, we have $\operatorname{cone}\theta \cong \mathfrak{t}_{\PP}I_P/\operatorname{soc}\mathfrak{t}_{\PP}I_P$. Applying the octahedral axiom to $\theta\circ\beta[d-1]=\real_{\PP}(\delta_P)[-1]$, we obtain the following commutative diagram of triangles:
\[\xymatrix{
&(\mathfrak{f}_{\PP}P)[d-1]\ar@{=}[r]\ar[d]_{\beta[d-1]}&(\mathfrak{f}_{\PP}P)[d-1]\ar[d]^{\real_{\PP}(\delta_P)[-1]}\\
(\mathfrak{t}_{\PP}I_P/\operatorname{soc}\mathfrak{t}_{\PP}I_P)[-1]\ar@{=}[d]\ar[r]& S[d-1]\ar[r]^{\theta}\ar[d]&\mathfrak{t}_{\PP}I_P\ar[r]\ar[d]&\mathfrak{t}_{\PP}I_P/\operatorname{soc}\mathfrak{t}_{\PP}I_P\ar@{=}[d]\\
(\mathfrak{t}_{\PP}I_P/\operatorname{soc}\mathfrak{t}_{\PP}I_P)[-1]\ar[r]&(\operatorname{rad}\mathfrak{f}_{\PP}P)[d]\ar[r]\ar[d]&\real_{\PP}(E_P)\ar[r]\ar[d]&\mathfrak{t}_{\PP}I_P/\operatorname{soc}\mathfrak{t}_{\PP}I_P\\
&(\mathfrak{f}_{\PP}P)[d]\ar@{=}[r]&(\mathfrak{f}_{\PP}P)[d]
}\]
Thus, by Theorem~\ref{thm:real-d-heart}, applying $\Theta_{\PP}$ to the third row yields the required triangle, where 
\[\Theta_{\PP}((\operatorname{rad}\mathfrak{f}_{\PP}P)[d])\in\X(\PP)\text{ and }\Theta_{\PP}(\mathfrak{t}_{\PP}I_P/\operatorname{soc}\mathfrak{t}_{\PP}I_P)\in\Y(\PP).\]
\end{proof}

%==================================================%
\section{Separating and splitting silting objects}\label{sec:sep-spli}
%==================================================%

We generalize the notions of separating and splitting 2-term silting complexes to the $(d+1)$-term setting.

An $s$-torsion pair $(\T,\F)$ in the $d$-extended module category of a proper $d$-truncated non-positive dg algebra is called \emph{split} if $\Hom(\F,\T[1])=0$. This is also equivalent to every indecomposable object in the extended module category belonging to $\T$ or $\F$.

\begin{definition}
  Let $A$ be a proper $d$-truncated non-positive dg algebra, and let $\PP$ be a $(d+1)$-term silting object in $\per A$. Let $B=\eend{d}{A}(\PP)$ be the $d$-truncated dg endomorphism algebra of $\PP$. Then
  \begin{enumerate}
    \item $\PP$ is called \emph{separating} if the $s$-torsion pair $(\T(\PP),\F(\PP))$ in $\emod{d}{A}$ is split, and
    \item $\PP$ is called \emph{splitting} if the $s$-torsion pair $(\X(\PP),\Y(\PP))$ in $\emod{d}{B}$ is split.
  \end{enumerate}
\end{definition}

The following lemma is useful.

\begin{lemma}\label{lem:proj-sim}
Let $A$ be a proper $d$-truncated non-positive dg algebra, and let $\PP$ be a $(d+1)$-term silting object in $\per A$. There are triangles in $\per A$:
\begin{equation}\label{eq:approx-tri}
 X_{r+1}\xrightarrow{g_{r+1}}A_r\xrightarrow{f_r}X_r\to X_{r+1}[1],
 \quad0\leq r\leq d,
\end{equation}
where $X_0=\PP$, $X_{d+1}=0$, $A_r\in\add A$ and
each morphism $f_r$ is a minimal right $\add A$-approximation.
Moreover, for every simple right $H^0(A)$-module $S$,
regarded as a dg $A$-module via the canonical projection $A\to H^0(A)$,
there are isomorphisms
\begin{equation}\label{eq:proj-sim}
 \Hom_{D(A)}(\PP,S[r])\cong\Hom_{D(A)}(A_r,S),
 \quad 0\leq r\leq d.
\end{equation}
\end{lemma}

\begin{proof}
The proofs of \cite[Lemmas~8.1, 8.2 and Proposition~8.3]{GZ} apply after interchanging the roles of $A$ and $\PP$ and replacing minimal left approximations by minimal right approximations.
\end{proof}

Separating silting objects can be shown to satisfy the equivalent conditions in Theorem~\ref{thm:equiv}.

\begin{proposition}\label{prop:separating}
  Let $A$ be a proper $d$-truncated non-positive dg algebra that is not acyclic, and let $\PP$ be a separating $(d+1)$-term silting object in $\per A$. Then
  \[\Hom_{\per A}(\PP,\PP[i])=0\quad\text{for all }i<-a,\]
  where $a$ is the largest non-negative integer such that $H^{-a}(A)\neq 0$. In particular, the equivalent conditions in Theorem~\ref{thm:equiv} hold.
\end{proposition}

\begin{proof}
  For any $1\leq r\leq d$, let $l_r=\max\{0,r-a-1\}$. Then
  \[0\leq l_r\leq d-1-a\text{ and }1\leq r-l_r\leq a+1.\]
  It follows that $A_r[l_r]\in D_{\mathrm{fd}}(A)^{[-a-l_r,-l_r]}\subseteq \emod{d}{A}$ for all $r$. We claim that $A_r[l_r]\in\F(\PP)$ for all $r$. Indeed, let $E$ be an indecomposable summand of $A_r$. By \eqref{eq:proj-equiv}, $H^0(E)$ is an indecomposable projective $H^0(A)$-module. Let $S$ be the simple top of $H^0(E)$. We regard $S$ as a dg $A$-module via the canonical morphism $A\to H^0(A)$. Then there exists a non-zero morphism from $E$ to $S$ in $D_{\mathrm{fd}}(A)$, which implies $\Hom_{D_{\mathrm{fd}}(A)}(A_r,S)\neq 0$. Note that both $E[l_r]$ and $S[l_r]$ are in $\emod{d}{A}$. Hence, by \eqref{eq:proj-sim} in Lemma~\ref{lem:proj-sim}, we have
  \[\Hom_{D_{\mathrm{fd}}(A)}(\PP,S[l_r][r-l_r])=\Hom_{D_{\mathrm{fd}}(A)}(\PP,S[r])\neq 0.\]
  Since $r-l_r\geq 1$, this implies $S[l_r]\notin\T(\PP)$. Since $(\T(\PP),\F(\PP))$ is split, we have $S[l_r]\in\F(\PP)$. Since $\Hom_{D_{\mathrm{fd}}(A)}(E[l_r],S[l_r])\cong\Hom_{D_{\mathrm{fd}}(A)}(E,S)\neq 0$, we have $E[l_r]\notin\T(\PP)$. Since $E[l_r]$ is indecomposable, we have $E[l_r]\in\F(\PP)$, which proves the claim.

  Fix a negative integer $i<-a$. For any $1\leq r\leq d$, since $A_r[l_r]\in\F(\PP)$ and $r+i-l_r\leq a+1+i\leq 0$, we have
  \[\Hom_{D_{\mathrm{fd}}(A)}(\PP,A_r[r+i])=\Hom_{D_{\mathrm{fd}}(A)}(\PP,A_r[l_r][r+i-l_r])=0.\]
  For the case $r=0$, since $A\in D_{\mathrm{fd}}(A)^{[-a,0]}$, we have $A_0[i]\in D_{\mathrm{fd}}(A)^{\geq 1}$. Since $\PP\in D_{\mathrm{fd}}(A)^{\leq 0}$, we have $\Hom_{D_{\mathrm{fd}}(A)}(\PP,A_0[i])=0$. Therefore, we have $\Hom_{D_{\mathrm{fd}}(A)}(\PP,A_r[r+i])=0$ for all $0\leq r\leq d$.

  Applying the functor $\Hom_{D_{\mathrm{fd}}(A)}(\PP,-)$ to the triangles in \eqref{eq:approx-tri} yields an exact sequence
  \[0=\Hom_{D_{\mathrm{fd}}(A)}(\PP,A_r[r+i])\to\Hom_{D_{\mathrm{fd}}(A)}(\PP,X_r[r+i])\to\Hom_{D_{\mathrm{fd}}(A)}(\PP,X_{r+1}[r+i+1]).\]
  Starting with $X_{d+1}=0$, descending induction gives $\Hom_{D_{\mathrm{fd}}(A)}(\PP,X_r[r+i])=0$ for all $0\leq r\leq d$. Since $X_0=\PP$, this completes the proof.
\end{proof}

An immediate corollary of the above result is the following, which generalizes \cite[Proposition~5.7]{BZ}.

\begin{corollary}\label{cor:separating}
   If $A$ is an ordinary finite-dimensional algebra, then any separating $(d+1)$-term silting object in $\per A$ is tilting.
\end{corollary}

For splitting silting objects, we have the following criterion;
compare \cite[Lemma~5.5]{BZ}.

\begin{lemma}\label{lem:splitting}
Let $A$ be a proper $d$-truncated non-positive dg algebra, and let $\PP$ be a $(d+1)$-term silting object in $\per A$. Then $\PP$ is splitting if and only if
\[
 \Hom_{D(A)}(\T(\PP),\F(\PP)[d+1])=0.
\]
\end{lemma}

\begin{proof}
This follows from the definition and the isomorphism~\eqref{eq:realization-Hom} in Theorem~\ref{thm:real-d-heart} for $i=1$.
\end{proof}

If $A$ is quasi-isomorphic to an ordinary finite-dimensional hereditary algebra, then we have $\Hom_{D(A)}(\emod{d}{A},\emod{d}{A}[d+1])=0$. Hence, in this case, by Lemma~\ref{lem:splitting}, every $(d+1)$-term silting object is splitting.

Since both $\T(\PP)$ and $\F(\PP)$ are functorially finite, by \cite[Theorem~5.14]{INP}, they also have almost split triangles.

\begin{proposition}\label{prop:AR-splitting}
Assume that $\PP$ is splitting. Then every almost split triangle in $\emod{d}{B}$ has its first three terms in either $\X(\PP)$ or $\Y(\PP)$, or else it is of the form
  \[\Theta_{\PP}(\mathfrak{t}_{\PP}I_P)\to\Theta_{\PP}((\operatorname{rad}\mathfrak{f}_{\PP}P)[d])\oplus\Theta_{\PP}(\mathfrak{t}_{\PP}I_P/\operatorname{soc}\mathfrak{t}_{\PP}I_P)\to\Theta_{\PP}((\mathfrak{f}_{\PP}P)[d])
    \to\Theta_{\PP}(\mathfrak{t}_{\PP}I_P)[1],
  \]
  where $P$ is an indecomposable object in $\add A$ with $P,P[d]\notin\add\PP$. Moreover, the functor $\Theta_{\PP}$ (resp. $\Theta_{\PP}(-[d])$) identifies the almost split triangles in $\T(\PP)$ (resp. $\F(\PP)$) with those in $\emod{d}{B}$ whose first three terms lie in $\Y(\PP)$ (resp. $\X(\PP)$).
\end{proposition}

\begin{proof}
  Since the pair $(\X(\PP),\Y(\PP))$ is split, every indecomposable object in $\emod{d}{B}$ belongs to $\X(\PP)$ or $\Y(\PP)$.

  The first statement follows from Lemma~\ref{lem:connecting2} and Proposition~\ref{prop:connecting}. To show the second statement, by Theorem~\ref{thm:real-d-heart}, it suffices to prove that any almost split triangle in $\X(\PP)$ (resp. $\Y(\PP)$) is also an almost split triangle in $\emod{d}{B}$. We only show this for $\X(\PP)$, since the case of $\Y(\PP)$ is dual. Let $X\xrightarrow{\alpha} Y\to Z\to X[1]$ be an almost split triangle in $\X(\PP)$. Let $f\colon X\to M$ be a non-section morphism in $\emod{d}{B}$. Write $M=M_1\oplus M_2$ with $M_1\in\X(\PP)$ and $M_2\in \Y(\PP)$. Then $f$ has two components $f_1\colon X\to M_1$ and $f_2\colon X\to M_2$, where $f_1$ is not a section and $f_2$ has to be zero since $\Hom_{D(B)}(\X(\PP),\Y(\PP))=0$. Thus, $f$ factors through $\alpha$. Since $Z$ is indecomposable, by \cite[Theorem~2.9]{INP}, $X\xrightarrow{\alpha} Y\to Z\to X[1]$ is an almost split triangle in $\emod{d}{B}$.
\end{proof}

\bigskip

\noindent \textbf{Use of AI:}
The authors found a proof of Proposition~\ref{prop:separating} with AI assistance and wrote it up themselves. AI tools were also used to correct grammar and spelling.

\end{document}